\documentclass[hidelinks,onefignum,onetabnum]{siamart220329}

\def\Ff{\mathcal{F}}

\def\P{\mathbb{P}}

\def\R{\mathbb{R}}

\def\COMMA{\,,}             
\def\PERIOD{\,.}            

\def\BARIT#1{{\bar {#1}}}

\def\Lapl#1{{\hat {#1}}}
\def\massat{\mathbb{M}}
\def\masscg{\bar{\mathbb{M}}}

\def\argmin{\textrm{ argmin}}

 \def\tr{[0,t_{f}]}
\def\bp{\mathbf{P}}

\def\ft{\bar{\zeta}_{[0,t_{f}]}}
\def\fe{\bar{\zeta}_{eq}}

\def\bx{Q}
 \def\bp{P}

\def\0o{\mathbb{O}}

\def\BARIT#1{{\bar {#1}}}

\def\BARV{\BARIT{U}}

\def\IP{U}
\def\PMF{{U}^{\mathrm{PMF}}}
\def\FPMF{F^{\mathrm{PMF}}}
\def\MVPMF{\mu^{\mathrm{PMF}}}

\def\BARMV{\BARIT{\mu}}
\def\BARF{\BARIT{F}}

\def\COP{\mathbf{\Pi}}

\usepackage{lipsum}
\usepackage{amsfonts}
\usepackage{graphicx}
\usepackage{epstopdf}
\usepackage{algorithmic}
\usepackage{subcaption}

\usepackage[markup=underlined]{changes}

\definechangesauthor[color=BrickRed]{GB}
\definechangesauthor[color=Green]{EK}
\definechangesauthor[color=NavyBlue]{VH}

\usepackage{bm,bbm}
\usepackage{tabu}
\usepackage{float}
 \usepackage[normalem]{ulem} 
 \usepackage{url}
\usepackage{color}
\usepackage{hyperref}

\ifpdf
  \DeclareGraphicsExtensions{.eps,.pdf,.png,.jpg}
\else
  \DeclareGraphicsExtensions{.eps}
\fi

\newsiamremark{hypothesis}{Hypothesis}
\crefname{hypothesis}{Hypothesis}{Hypotheses}
\newsiamthm{claim}{Claim}

\headers{Bottom-up Transient Time Models}{G. Baxevani, V. Harmandaris, E. Kalligiannaki, and I. Tsantili}

\title{Bottom-up Transient Time Models in Coarse-graining Molecular Systems\footnote{The content of the manuscript, in whole or in part, is not submitted, accepted, or published elsewhere, including conference proceedings.}}

\author{Georgia Baxevani\thanks{Department of Mathematics and Applied Mathematics, University of Crete, Heraklion, GREECE 
  (\email{g.baxevani@iacm.forth.gr}, \email{harman@uoc.gr})}
  \and Vagelis Harmandaris\footnotemark[2] \footnotemark[5] \thanks{Computation-Based Science and Technology Research, Center, The Cyprus Institute, 2121 Nicosia, Cyprus}
  \and  Evangelia Kalligiannaki\footnotemark[5] \thanks{Department of Statistics and Actuarial-Financial Mathematics, University of the Aegean, GREECE (\email{ekalligiannaki@aegean.gr})}
  \and Ivi Tsantili\thanks{Institute of Applied and Computational Mathematics FORTH, Heraklion, GREECE (\email{ivi.tsantili@iacm.forth.gr})  } 
  }
\usepackage{amsopn}

\begin{document}

 \maketitle

\begin{abstract}
 This work presents a systematic methodology for describing the transient dynamics of coarse-grained molecular systems inferred from all-atom simulated data. We suggest Langevin-type dynamics  where the coarse-grained interaction potential depends explicitly on time to efficiently approximate the transient coarse-grained dynamics. We apply the path-space force matching approach at the transient dynamics regime to learn the proposed model parameters. In particular, we parameterize the coarse-grained potential both with respect to the pair distance of the coarse-grained particles and the time, and we obtain an evolution model that is explicitly time-dependent. Moreover, we follow a data-driven approach to estimate the friction kernel, given by appropriate correlation functions directly from the underlying all-atom molecular dynamics simulations. To explore and validate the proposed methodology  we study a benchmark system of a moving particle in a box. We examine the suggested model's effectiveness in terms of the system's correlation time and find that the model can approximate well the transient time regime of the system, depending on the correlation time of the system. As a result, in the less correlated case, it can represent the dynamics for a longer time interval. We present an extensive study of our approach to a realistic high-dimensional water molecular system. Posing the water system  initially out of thermal equilibrium we collect trajectories of  all-atom data for the, empirically estimated, transient time regime. Then, we infer the suggested model and strengthen the model's validity by comparing it with simplified Markovian models.
\end{abstract}




\section{Introduction}
Simulation  of complex molecular systems is a very intense research area \cite{FrenkelSmitBook}. The main challenge in this field is to predict the properties of materials and to provide a direct quantitative link between chemical structure at the molecular level and measurable structural and dynamical quantities over a broad range of length and time scales. On the microscopic (atomistic) level, detailed all-atom (AA) molecular dynamics (MD) or Monte Carlo (MC) simulations allow the description of molecular systems at the most detailed level \cite{AllenTildesleyBook,FrenkelSmitBook},\cite{MulPlat2002,Soper1996,Harmandaris2003a,TheodorouBook}. However, these standard computational techniques encounter temporal and spatial limitations resulting in an unfeasible investigation of phenomena occurring at larger scales. 

Consequently, much attention has been paid to developing \textit{coarse-grained} (CG) models. On the mesoscopic level, CG models have proven to be very efficient means in order to increase the length and time scales accessible by simulations \cite{Soper1996,MulPlat2002,briels},
\cite{Harmandaris2006a,Harmandaris2009a,Johnston2013}.
These models reduce the computational cost by including only a relatively small number of the most essential degrees of freedom in the computation. 
Much effort has been devoted to the  approximation of  the many-body \textit{potential of mean force} (PMF). As a result, a  variety of methods have been developed, such as the Boltzmann inversion \cite{Soper1996,MulPlat2002,Reith2003}, force matching (FM) (or Multi-Scale CG) \cite{IzVoth2005a,Voth2008a,Voth2010, Noid2013, Andersen2012}, relative
entropy (RE) \cite{Shell2008,Katsoulakis2013,Zabaras2013,KHKP2015,Shell2016, Kalligiannaki2016e}, and inverse MC \cite{ LyubLaa2004,Lyubartsev2010, Lyubartsev2003}. Beyond these methods which develop pair-wise models aiming at matching the structure of the AA model, a different approach lies in the development of local density potentials to include the many-body terms. This methodology can serve as a basis to improve predictions on both structure and thermodynamics \cite{Moore2016,doxastakis2020}.

While many existing CG models have demonstrated their capability to recover structural properties, a systematic framework incorporating dynamic properties is still challenging. Therefore, there is a significant interest in studying CG models that resemble the dynamic behavior, both at equilibrium and out-of-equilibrium conditions \cite{VothDyn2006,VothDyn2015,DiffusionMapKevrekidis2008,Katsoulakis2013,HKKP2016,KalligiannakiUQ, VanderVegt2018, Schilling_2022}. The well-known \textit{Generalized Langevin equation} (GLE), resulting through the   $\textit{Mori-Zwanzig}$ projection operator formalism, describes the evolution of the CG particles at equilibrium conditions  \cite{Evans-Morris,Mori1965,Zwanzig1973}. 
However, due to the intractability of the GLE, which lies mainly on the evaluation of the memory kernel, there have been numerous attempts to approximate the CG dynamics using alternative approaches \cite{Li_Liu2014,Chorin,lelievre2010,Buscalioni2010,ChunGLE2016,tuckerman2011}. The most widely used models are Markovian approximations of the GLE, in which the memory kernel is approximated by a Dirac delta function and a \textit{Langevin equation} is recovered  \cite{DiffusionMapKevrekidis2008,Espanol2006,Coffey2012}. 
Furthermore, in recent studies \cite{Lei14183,ChunGLE2016,Baker2020,LIN2021,Kevrekidis2019} data-driven methods have been introduced using a variety of approaches to representing non-Markovian dynamics.  

Dynamics models, under the Markovian approximation, are also recovered with variational inference methods 
 \cite{Katsoulakis2013,HKKP2016}. In these works, the authors address a \textit{path-space} variational inference problem in which they derive an approximation of the exact CG model by optimizing the information content with respect to time-series data.
Recovering the correct  dynamics relies on the availability of large datasets corresponding to the equilibrium dynamics regime \cite{KalligiannakiUQ}. 
However, for complex molecular systems, there is often an extremely large computational cost in generating a large dataset of atomistic configurations. For example, polymer-based materials are characterized by long relaxation times, and thus sampling a large number of atomistic configurations is not feasible \cite{Harmandaris2003a,Harmandaris2009a}.

In this context, two challenging directions associated with still-open questions arise. First, how can we estimate the dynamics of high-dimensional molecular systems for which it is not feasible to perform extensive atomistic simulations that allow deriving equilibrium data? Second, how can we study the structural and dynamic   behavior of molecular systems under out-of-equilibrium conditions? 
In both cases, there is no available information corresponding to the invariant (or the steady state) measure; thus, in principle, we can not efficiently estimate the PMF and the corresponding friction forces.

The goal of the current study is to approximate the evolution, with respect to both the structure and dynamics, of CG model systems when they are at the \textit{transient regime} and thus have not reached equilibrium. A system is said to be in a transient regime when, after disturbance, it moves from one steady state to another or has not yet reached a steady state. 
We are interested in how the CG system, which initially, at $t=0$, is out of equilibrium, evolves in the time interval $\tr$ corresponding to the transient regime. Thus, we propose novel stochastic dynamical equations to approximate the GLE at the transient time regime.

The study of the behavior of stochastic systems at transient time regimes is a very active field. For example, recent works refer to a detailed sensitivity analysis for transient stochastic systems \cite{Arampatzis2015} and to model reduction and transient dynamics \cite{Babaee,zhang2011reconstructing}. In the latter work, the authors present a variational principle and develop a reduced description of the transient dynamics. Inference methods for far from equilibrium systems with Langevin dynamics are also studied lately in \cite{nature_Genkin_2021}.
 
In this work, we go beyond the widely used Langevin models as the standard Markovian models fail to investigate the behavior of "bottom-up" CG models at the transient time regimes. 
To address the above challenges, i.e. the short-time observations and the inherent memory, we approximate the short-time non-Markovian dynamics by Markovian dynamics with a time-dependent force field. Specifically, we approximate the GLE, which has memory terms, by Langevin dynamics with a time-dependent force field that could incorporate the memory into the CG dynamics. 
Our approach is based on a recent work, where we proved that one can derive a time-depended potential valid for the time interval of observation and verified numerically that it reduces to the PMF as simulation time increases \cite{BAXEVANI201959}.

We parametrize the CG potential with respect to \textbf{(a)} the pair distance between the CG particles and (\textbf{b)} the simulation time. To learn the proposed model parameters, we apply the \textit{path-space force-matching} (PSFM) approach at the transient dynamics regime \cite{HKKP2016}. We study the effect of the actual time-dependent force field, as well as the corresponding friction forces, on the system evolution at the transient regime by probing directly both its structure and dynamics and comparing directly with the underlying data of the microscopic (atomistic) simulations.
The proposed approach is applied in a benchmark example concerning the dynamics of one particle in a box with time-dependent potential, as well as in a more realistic model of liquid water under non-equilibrium conditions.  

We should note that dynamics with a time-dependent interaction potential were introduced recently in the case of non-stationary dynamics to describe CG dynamics. Specifically,  Meyer et al. \cite{doi:10.1063/1.5006980}, and Vrugt et al. \cite{PhysRevE.99.062118} employed the formalism of time-dependent projection operator techniques and derived an equation of motion of the CG particles for the non-stationary case, the non-stationary GLE. It resembles the GLE but differs from the GLE in the explicit dependence of the CG interaction potential, the memory kernel, and the fluctuation force on time. In contrast, our work concerns estimating the transient CG dynamics for stationary processes, that is assuming that as time increases the invariant probability measure is ensured, and the time-dependent interaction potential approaches the PMF.


The remainder of this paper is organized as follows. In section \ref{AtCG}, we describe the microscopic and CG models and the corresponding dynamics. Section \ref{proposed transient model} describes the proposed CG model for describing the system's transient dynamics. The data-driven methods for the approximation of the time-depended CG interaction potential and the friction coefficient are presented in \ref{inferring cg dynamics}. Next, in section \ref{benchmark}, we examine a benchmark application;
a GLE describing the movement of one particle in a box approximated by the Langevin equation with time-dependent potential.
The molecular models and details about the atomistic and CG simulations and the results for a high-dimensional liquid water system are given in section \ref{results}. Finally, in the last section \ref{discussion-conclusions}, we summarize and discuss the main results of the present work.

\section{Atomistic and coarse-grained systems}\label{AtCG}
Assume a prototypical problem of $N$ (classical) particles in a box of volume $V$ at temperature $T$.
Let $q=(q_1,\dots,q_N) \in \R^{3N}$ and $p=(p_1,\dots,p_N) \in \R^{3N}$ describe the position and the momentum vector, respectively, of the $N$ particles in the atomistic (microscopic)  description, with potential  energy $\IP(q)$. We consider the evolution of the $N$ particles described by the process $\{q_{t},p_{t}\}_{t\ge 0}\in \mathbb{R}^{6N}$, given the initial conditions $\{q_{0},p_{0}\}\in \mathbb{R}^{6N}$, satisfying
\begin{equation}
\begin{cases}
    dq_t=\massat^{-1}p_tdt,\\
    dp_t=f(q_t)dt-\zeta \massat^{-1} p_tdt+\sigma dW_t, \quad t>0,\\
      (q_{0},p_{0}) \sim \nu_0(dq,dp)\COMMA
  \end{cases}
 \end{equation} 
a Hamiltonian system coupled with a thermostat. We denote  $f(q)$ the force vector applied to particles such that $f=(f_{1},\dots,f_{N})\in \R^{3N} $.
$\massat=diag[m_{1}I_{3},..,m_{N}I_{3}]$ is the mass matrix, with $m_{i},\: i=1,..,N$ the mass of the atom $i$ and $I_{3}$  the 3-dimensional identity matrix. $W_t$ denotes the $3N-$dimensional Wiener process. $\zeta \in \mathbb{R}^{3N\times 3N}$ and $\sigma \in \mathbb{R}^{3N\times 3N}$ are the friction and the diffusion coefficients, respectively. The notation $(q_{0},p_{0}) \sim \nu_0(dq,dp)$ means that the initial variables $(q_{0},p_{0})$ are distributed according to the probability measure $\nu_0$. 

The diffusion and friction coefficients satisfy the fluctuation-dissipation theorem (FDT) $\sigma\sigma^{tr}=2\beta^{-1}\zeta$,  where $\beta=\frac{1}{k_b T}$ and $k_b$ is the Boltzmann constant and $\cdot^{tr}$ denotes matrix transpose. The FDT expresses the balance between friction and noise in the system. This balance is crucial to ensure that the system reaches a thermal equilibrium state at long times as well as to obtain the Gibbs canonical measure at equilibrium, given by
\begin{equation*}
    \mu(dq)=Z^{-1}e^{-\beta U(q)}dq\COMMA
\end{equation*}
where $Z= \int_{\R^{3N}}e^{-\beta \IP(q)} dq$ is  the partition function. 
 
\textit{Coarse-graining} is a standard methodology to overcome the difficulties that lie in the fully detailed AA simulations and lead to
temporal and spatial limitations. CG models reduce computational costs by including only a relatively small number of the most important degrees of freedom. Thus, CG models allow access to processes on a broader range of spatiotemporal scales by providing a lower-resolution description of a given molecular system. Specifically, coarse-graining is the application of a mapping (CG mapping)
$$\COP:\R^{3N} \to \R^{3M}$$
\begin{equation*}
      q  \mapsto \COP (q)\in \R^{3M}
\end{equation*}
on the microscopic state space, determining the $M(<N)$ CG particles as a function of the atomistic configuration $q$. We denote by $ \bx = (\bx_1,\dots,\bx_M)\in\mathbb{R}^{3M}$ any point in the CG configuration space, and use the bar $"\ \BARIT{\ }\ "$ notation for quantities on the CG space.  

The probability that the CG system has configurations in $d\bx$ is given by 
\begin{equation*}
   \MVPMF(dQ) =  \int_{\Omega(\bx)}\mu(dq)\,dQ=Z^{-1}\int_{\Omega(\bx)}e^{-\beta U(q)}dq \, dQ\COMMA
\end{equation*}
where $\Omega(\bx)=\{q\in \mathbb{R}^{3N}:\mathbf{\Pi}(q)=\bx\}$.
The  $M-$body PMF is thus
\begin{equation}
\PMF(\bx)=-\frac{1}{\beta}\ln\int_{\Omega(\bx)}e^{-\beta U(q)}dq. \label{analytical form of PMF}
\end{equation}
 
The PMF calculation is a task as complex and costly as calculating expectations on the microscopic space. Thus, one seeks a practical potential function $\BARV(\bx)$ that approximates as well as possible the PMF. Estimating a potential $\BARV(\bx)$ is the goal of the inference methods for systems at equilibrium, i.e., structural-based methods, (FM, RE minimization, inverse MC) \cite{MulPlat2002,briels,IzVoth2005a,Shell2008,Shell2009,IzVoth2005,KHKP2015,Kalligiannaki2016e,Lyubartsev2003,Lyubartsev2010,Voth2008b,Harmandaris2006a,IzVoth2005a,Voth2008a,Voth2010,Noid2011,Noid2013,Zabaras2013,Shell2016,Shell2011,Andersen2012}. In all these methods, one proposes a family of potential functions  in a parameterized, $ \BARV( \bx;\phi)  $, $\phi\in  \Phi$, or in a functional  form $  \BARV( \bx)$ and seeks for the optimal $\BARV^*(\bx)$ that best approximates the PMF.
We denote by 
\begin{equation*}\label{CGmeasure}
\BARMV(d\bx) =\ \BARIT Z^{-1}  e^{ -\beta \BARV( \bx)}d\bx\COMMA
\end{equation*}
the  equilibrium probability measure at the CG configurational space 
where $\BARIT Z =  \int e^{ -\beta \BARV( \bx)} d\bx$. The knowledge of an effective CG interaction potential $ \BARV(\bx)$ is enough to approximate  structural properties efficiently for systems at equilibrium. However, dynamic properties, e.g., diffusion or, more generally, time correlations are, in principle, not recovered correctly from MD simulations using only the CG potential, i.e., only the conservative force. That is, the CG procedure eliminates the degrees of freedom that should appear in the CG dynamics. Thus, we need to estimate the friction kernel function and the random forces arising in \eqref{eq:GLE} to find the coarse particles' dynamics.

The real dynamics of the coarse particles are described by the process \\\mbox{$\{Q_{t},P_t\}_{t\geq0}\in\mathbb{R}^{3M}$} satisfying the GLE
\begin{equation}\label{eq:GLE}
    \begin{cases}
    d\bx_t = \masscg^{-1} \bp_t dt \\
    d\bp_t = \FPMF(\bx_t)dt - \int_{0}^{t}\bar{\zeta}(t-\tau)\masscg^{-1} \bp_\tau d\tau + F_R(t)
    \COMMA \quad t>0 \COMMA\\
     (\bx_{0},\bp_{0})\sim\bar{\nu}_0(d\bx,d\bp)\COMMA
    \end{cases}
\end{equation}
which is derived via a projection operator technique, the   \textit{Mori-Zwanzig} theory \cite{tuckerman2011,Chorin}. $\FPMF(\bx) =-\nabla \PMF(\bx)$ is the conservative force 
and $\bar{\zeta}(t)$ is the memory  kernel for the friction force. The random force $F_R(t)$ denotes fluctuating forces from the variables eliminated through the CG process. $\masscg$ is the mass matrix of the CG particles. 
 If the FDT is satisfied for the coarse space dynamics, 
\begin{equation*}
 \langle F_R(t)\cdot  F_R(t') \rangle  = 2 \beta^{-1}\bar{\zeta}(t-t')\COMMA\: t>t'
\end{equation*}
  a thermal equilibrium state is reached and the CG configurations are distributed according to $\MVPMF(d\bx)$, when $t\to \infty$ or when $\bar{\nu}_0(d\bx,d\bp) \propto \MVPMF(d\bx)  e^{-\frac{\beta}{2}\bp \masscg^{-1} \bp^{tr}}   d\bp$.

Although the CG dynamics is known given by \eqref{eq:GLE}, its practical use is quite challenging \cite{Lei14183,ChunGLE2016}. In principle, we do not know the closed form of the force and the friction kernel; thus, approximations are necessary. Moreover, the memory term's direct evaluation is costly because it requires the history of the CG variables at every time step and the associated numerical integration. To overcome these limitations simplified approximate models are proposed. Next, we approximate the actual dynamical process $\{Q_{t},P_t\}_{t\geq0}\in\mathbb{R}^{3M}$ with the Markovian process $\{\hat{Q}_{t},\hat{P}_t\}_{t\geq0}\in\mathbb{R}^{3M}$ which satisfies the Langevin equation, 
\begin{equation}\label{eq:Langevin equation}
    \begin{cases}
    d\hat{\bx}_t = \masscg^{-1} \hat{\bp}_t dt \\
    d\hat{\bp}_t = \BARF(\hat{\bx}_{t})dt -\bar{\zeta}_{0}\masscg^{-1} \hat{\bp}_{t}dt + \bar{\sigma}_{0} d\bar{W}_t\COMMA \quad t>0 \COMMA\\
     (\hat{\bx}_{0},\hat{\bp}_{0})\sim\bar{\nu}_0(d\hat{\bx},d\hat{\bp})\COMMA
    \end{cases}
\end{equation}
that is a Markovian approximation of the GLE. This is the case for $\langle F_R(t)\cdot F_R(t') \rangle = 2 \beta^{-1}\bar{\zeta}_{0}\ \delta(t-t') $, where $\delta$ denotes the Dirac delta function.
The friction coefficient of the CG system is a constant, $\bar{\zeta}_{0}$ related to the diffusion coefficient through the FDT $\bar{\sigma}_{0}\bar{\sigma}_{0}^{tr}=2\beta^{-1}\bar{\zeta}_{0}$. The force $\BARF(Q)$ is an approximation of the $\FPMF(Q)$,  chosen as $\BARF(Q) = -\nabla \BARV(Q)$.


In this work, we go beyond the above Langevin models to encounter two classes of problems that often appear in the  CG modeling of molecular systems. First, for the case when the system is at a transient state and has not reached equilibrium yet. In such a case, no available information corresponds to the invariant measure; thus, we can not estimate the $\FPMF$. The second class of problems we address is when the Markovian models given by \eqref{eq:Langevin equation} fail to capture the effect of the actual friction forces on the system dynamics.

Thus, to model CG systems for such situations we propose Langevin type dynamics with time-dependent drift coefficients. 

\section{Transient time dynamics models} \label{proposed transient model}
We are interested in finding the CG dynamics when the system is at a transient time regime. Thus, we assume next that we observe the system in a short-time interval $\tr$. 
Then,  we look for an approximation of the GLE \eqref{eq:GLE} valid at the transient regime in the form of Markovian Langevin-type dynamics, \eqref{eq:Langevin equation}. 

The novelty in the new proposed dynamics lies in the fact that the force $\bar{F}$   depends explicitly on time. That is, we propose Langevin-type dynamics characterized by the time-depended force $\bar{F}(\tilde{\bx},t)$ and the constant friction coefficient $\bar{\zeta}_{0}$, defining the process $\{\tilde{Q}_{t},\tilde{P}_t\}_{t\geq0}\in\mathbb{R}^{3M}$ by
\begin{equation}
    \begin{cases}
    d\tilde{\bx}_t = \masscg^{-1}\tilde{\bp}_t dt \\
    d\tilde{\bp}_t = \bar{F}(\tilde{\bx}_{t},t)dt -\bar{\zeta}_{0}\masscg^{-1} \tilde{\bp}_{t}dt + \bar{\sigma}_{0} d\bar{W}_t\COMMA \quad t>0 \COMMA\\ 
    (\tilde{\bx}_{0},\tilde{\bp}_{0}) \sim \bar{\nu}_0(d\tilde{\bx},d\tilde{\bp})\PERIOD
    \end{cases} \label{eq:Langevin_time_d}
\end{equation}
   
We assume that the system's initial state is out of equilibrium. Such a model aims at approximating the memory part of the friction forces with time-varying interacting forces while preserving the actual dynamics of the system for the time interval  $\tr$. 
To estimate the force $\bar{F}( \tilde{\bx},t) $ and friction coefficient $\bar{\zeta}_{0}$ we apply   data-driven approaches, which we  describe in detail in section \ref{inferring cg dynamics}.
 
 The above model \eqref{eq:Langevin_time_d} belongs to a broader class of stochastic differential equations (SDEs), for the stochastic process $\{{X}_t\}_{t\ge 0}\in \R^d$,
\begin{equation}
dX_t={b}(t,X_{t})dt+{\sigma}(t,X_{t})dB_{t}, \quad t>0, \quad X_0= x_0\COMMA
    \label{eq:general_SDE}
\end{equation}
where 
$ {b}:[0,T]\times\R^{d}\to\R^{d}$   ${\sigma}:[0,T]\times\R^{d}\to\R^{d}\times\R^{n}$ are measurable vector-valued functions, $B_{t}$ is the standard Brownian motion in $\R^{n}$. The initial condition $x_0$ is a random variable independent of the Brownian motion $B_{t}$.

Equation \eqref{eq:Langevin_time_d} is a special case of equation \eqref{eq:general_SDE} for the process $  X_{t}=(\tilde{\bx}_{t},\tilde{\bp}_{t}):[0,T]\to\R^{d}$ for $d=6M$, 
${b}(t,x)=(\masscg^{-1} \bp ,\bar{F}( \bx ,t) -\bar{\zeta}_{0}\masscg^{-1}  \bp )$, and ${\sigma}(t,x) = \begin{bmatrix} 0&0\\ 0&\bar{\sigma}_{0} \end{bmatrix}$ where $x=(\bx,\bp)$. 
Equation \eqref{eq:general_SDE} admits a unique strong solution, \cite{pavliotis2008,oksendal2003}, $X_{t}:[0,T]\to\R^{d}$, defined on the probability space $(\Omega, \Ff, \P)$ when  ${b}(t,x),{\sigma}(t,x)$ are globally Lipschitz vector and matrix fields and follow the linear growth condition, i.e., there exists a positive constant $C$ such that for all ${x},y$ in $\R^{d}$
\begin{equation*}
|{b}(t,x)-{b}(t,y)|+|{\sigma}(t,x)-{\sigma}(t,y)| \leq C(|x-y|)
    \label{Glob_L1}
\end{equation*} and
\begin{equation*}
|{b}(t,x)|+|{\sigma}(t,x)| \leq C(1+|x|).
    \label{Glob_L2}
\end{equation*}
Here $|{\sigma}(t,x)|^2 = \sum_{ij} |\sigma_{ij}(t,x)|^2$. 
Therefore to ensure existence and uniqueness   for our model dynamics \eqref{eq:Langevin_time_d},
we   require that $|\bar{F}( \bx,t)|$  has linear growth and is Lipschitz in $\bx$. That is  since the diffusion coefficient is constant  and the drift is linear in $P$.

\section{Inferring coarse dynamics}\label{inferring cg dynamics}
We seek a CG model that best represents the atomistic reference system at a transient regime, ideally best in terms of its structural and dynamic properties. We follow data-driven approaches to estimate the CG interactions and the friction coefficient in the proposed CG model given in \eqref{eq:Langevin_time_d}.

In particular, for estimating the CG interactions we follow the PSFM approach at transient dynamics regimes. We consider that a stochastic process describes the evolution of the CG particles, see subsection \ref{PSFM subsection}. The PSFM is a generalization of the FM method \cite{IzVoth2005a,Noid2011,IzVoth2005} developed in previous works \cite{HKKP2016,KALLIGIANNAKI2018} extending the FM method for systems at equilibrium to non-equilibrium systems. It has been proved that when the CG forces are time-independent and the system  has reached equilibrium, the FM and the PSFM methods are equivalent, i.e., they produce the same effective CG interaction potential \cite{HKKP2016}. To estimate the friction kernel  we  also follow a data-driven approach based on the Markovian assumption,  calculating  appropriate correlation functions directly from the underlying AA MD simulations \cite{VothDyn2006}, see subsection \ref{representation of the friction kernel subsection}.

\subsection{Path-space force matching method} \label{PSFM subsection}
The PSFM method is entirely data-driven and it  involves information along the  path space. This path-space approach concerns the dynamic equilibrium, transient, and non-equilibrium regime. In this method, the atomistic data, which are entire paths, are correlated, i.e., they arise from observing the system over a continuous period of time.

Observing the system in a short time interval of the transient regime, $\tr$, the optimization problem is to find the optimal $\bar{F}^{*}(\bx)$, $\bar{F}^{*}:\mathbb{R}^{3M}\to\mathbb{R}^{3M}$ such that the mean square error,

\begin{equation}
\chi^{2}_{\tr}=\mathbb{E}_{\mathcal{P}_{\tr}}\left[\int_{0}^{t_{f}} \|\COP f(q_{s})-\bar{F}(\COP q_{s})\|^{2}ds\right]\COMMA \label{FM ps}
\end{equation}
is minimized \cite{HKKP2016, KHKP2015}. In other words, we minimize the average difference between the mapped atomistic forces $F(q)=\COP f(q)$  and the corresponding proposed CG forces $\bar{F}(\COP q)$, where $\|\cdot \|$ denotes the Euclidean norm in $\mathbb{R}^{3M}$ and $\mathbb{E}_{\mathcal{P}_{\tr}}[\cdot]$ averages with respect to the path-space probability measure $\mathcal{P}_{\tr}$ associated with the exact process $\left\{ q_{t},p_{t}  \right\}_{0\leq t\leq t_{f}}$. As $F(q)=\COP f(q)$, we denote the vector of forces exerted on the CG particles, which are given as a function of the microscopic forces. 
If we consider the CG mapping as the center of mass of a group of atoms, the force exerted on the $I-$th CG particle is given by $F_{I}(q)=\sum_{j \in I  }f_{j}(q)$, $I=1,...,M,$ where $f_{j}(q)$ is the atomistic force acting on atom $j$.
Hence, the PSFM is seeking the optimal force $\bar{F}^{*}(\bx)$ as the minimizer of the mean square error $\chi^{2}_{\tr}$ over a set of proposed CG forces $\bar{F}(\bx) = -\nabla  \bar{U}(\bx) $ defined on a suitable function space.

In principle, there are many-body contributions to the PMF, such as two-body, three-body, etc. Our study assumes a two-body effective pair potential to approximate the PMF
\begin{equation*}
  \bar{U}(\bx)=\sum_{I=1}^{M}\sum_{J>I}^{M} \bar{u}(r_{IJ})\COMMA
\end{equation*}
where $\bar{u}(r)$ is the potential that corresponds to the pair interaction between CG particles $I$ and $J$, at distance $r_{IJ}$, such that $r_{IJ}=|\bx_{I}-\bx_{J}|,\:I\neq J,\:I,J=1,..,M.$ 
When the pair forces depend only on the CG positions and are time-independent, we parameterize the CG forces with respect to pair distance $r$ between the CG particles, $\bar{F}(\bx;\phi),\:\phi\in \Phi$, and 
\begin{equation*}
    \bar{U}(\bx;\phi)=\sum_{I=1}^{M}\sum_{J>I}^{M}\bar{u}(r_{IJ};\phi)\PERIOD
\end{equation*}
As an example of parametrization, one can write the pair CG forces as a linear combination of non-parametric basis functions, such as linear or cubic splines \cite{Kalligiannaki2016e}.
In this case, the empirical analog of the minimization problem \eqref{FM ps}, which provides the optimal parameter set $\phi^{*}\in \Phi$, is 

 \begin{equation}\label{discrete PSFM}
      \phi^{*}_{\tr}=  \argmin_{\phi\in \Phi}\frac{1}{3M}\frac{1}{n_{t}}\frac{1}{n_{p}}\sum_{n=1}^{n_p}\sum_{i=1}^{n_{t}}\sum_{I=1}^{M}||F_{I}(q_{i,n})-  \bar{F}_{I}(\COP q_{i,n};\phi)||^{2}\COMMA \qquad \textbf{PSFM}\ 
 \end{equation}
where  $\bar{F}_{I}(\bx;\phi)=-\nabla_{\bx_I}\bar{U}(\bx;\phi)  $, and $\{q_{i,n}\},i=1,..,n_{t},\:n=1,..,n_{p}$, is the generated position data set which consists of $n_{p}$ independent (uncorrelated) paths with $n_{t}$ 
configurations per path, sampled by atomically detailed simulations on $\tr$. In \ref{discrete PSFM}, $||\cdot||$ denotes the Euclidean norm in $\R^{3}$. The optimal parameter set $\phi^{*}_{\tr}$ is   implicitly time-depended through the time-interval $\tr$ for which  we apply the PSFM method. It is noteworthy that if $t_{f}\to\infty$, due to the ergodicity theorem, even if $n_{p}=1$, we produce the PMF, i.e., at long times, the PSFM method reduces to the FM method, as has been recently verified \cite{BAXEVANI201959}.
 
A sub-case of this method is considered when we solve the PSFM method and derive a valid CG potential at a specific instantaneous time $t$. In this case, we use as a sample the configurations at the particular instantaneous time, $t\in \tr$, extracted from many different paths. The corresponding optimization problem is 
\begin{equation}\label{PSFM instant times}
       \phi^*_{t}=  \argmin_{\phi\in \Phi}\frac{1}{3M}\frac{1}{n_{p}}\sum_{n=1}^{n_p}\sum_{I=1}^{M}||F_{I}(q_{n}(t))-  \bar{F}_{I}(\COP q_{n}(t);\phi)||^{2}\COMMA 
 \end{equation}
where the number of configurations of each path is now $n_{t}=1$. In this way, we can produce the time evolution of the CG potential by solving this method for many instantaneous times of the transient short-time-interval $\tr$. Similarly, if $t\to\infty$, the system reaches equilibrium, and the derived CG potential approaches the PMF.
 
The current study considers that the CG interaction potential depends on time explicitly. Thus, we parameterize the time-dependent potential with respect to the pair distance $r$ of the CG particles and the time $t$ simultaneously. Therefore, the time-dependent parametrized CG potential is written in the form
\begin{equation*}
    \bar{U}(\bx,t;\theta)=\sum_{I=1}^{M}\sum_{J>I}^{M}\bar{u}(r_{IJ},t;\theta)\PERIOD
\end{equation*}
The main difference between the parameter sets $\phi$, and $\theta$ is that the parameter set $\phi$ is used in representing the pair CG forces (or the pair CG interaction potential) with respect to the pair distance $r$ of the CG particles. In contrast, we use $\theta$ to represent the same quantity with respect to $r$ and time $t$ simultaneously. 
Thus, to find an effective time-dependent CG potential, we solve the following optimization problem that provides the optimal parameter set $\theta^*$

  \begin{equation}\label{splines exp. on t and r, opt.problem}
    \begin{split}
  \theta^*_{\tr}= argmin_{\theta\in \Theta}\frac{1}{3M}\frac{1}{n_{t}}\frac{1}{n_{p}}\sum_{n=1}^{n_p}\sum_{i=1}^{n_{t}}\sum_{I=1}^{M}||F_{I}(q_{i,n})-  \bar{F}_{I}(\COP q_{i,n},t;\theta)||^{2}\PERIOD
    \end{split}
\end{equation}


Next, we present the parametric model of the CG force field. We expand the pair forces $\bar{f}(r,t;\theta)=-\frac{\partial \bar{u}(r,t;\theta)}{\partial r}$ on two-dimensional splines with respect to $t$ and $r$
\begin{equation*}
  \bar{f}(r,t;\theta)  =\sum_{d=1}^{N_{d}}\sum_{b=1}^{N_{b}}\theta_{db}\psi_{d}(r)\chi_{b}(t)\COMMA
\end{equation*}
where $N_{d}  N_{b}$ is the total number of basis functions, $\psi(r)$ and $\chi(t)$ for $r$ and $t$, respectively, and $\theta$ is the total number of parameters such that\\
$\theta=(\theta_{11},\theta_{12},\dots,\theta_{1N_{b}},\theta_{21},\dots,\theta_{2N_{b}},\dots,\theta_{N_{d}1}\dots\theta_{N_{d}N_{b}})$.
The CG force acting on CG particle $I$ is
\begin{align*}
    \bar{F}_{I}(\bx,t;\theta)=&\sum_{J>I}^{M}\sum_{d=1}^{N_{d}}\sum_{b=1}^{N_{b}}\theta_{db}\psi_{d}(r_{IJ})\chi_{b}(t)\nabla_{\bx_{I}}r_{IJ}\\=&\sum_{d=1}^{N_{d}}\sum_{b=1}^{N_{b}}\theta_{db}\chi_{b}(t)\cdot\sum_{J>I}^{M}\psi_{d}(r_{IJ})\nabla_{\bx_{I}}r_{IJ}\PERIOD
\end{align*}
Reordering the parameter indices 
\begin{align*}
    \bar{F}_{I}(\bx,t;\theta')=&\sum_{s=1}^{N_{d}\times N_{b}}\theta'_{s}\chi_{b}(t)\cdot \sum_{J>I}^{M}\psi_{d}(r_{IJ})\nabla_{\bx_{I}}r_{IJ}\\=&\sum_{s=1}^{N_{d}\times N_{b}}\theta'_{s}\cdot J_{I;s}(\bx)\COMMA
\end{align*}
where $\theta'=(\theta'_{1},\dots,\theta'_{N_{b}},\theta'_{N_{b}+1},\dots,\theta'_{2 N_{b}},\dots,\theta'_{N_{d} N_{b}})$  and \\ $J_{I;s}(\bx)=\chi_{b}(t)\cdot \sum_{J>I}^{M}\psi_{d}(r_{IJ})\nabla_{\bx_{I}}r_{IJ}$. 
Thus, the minimization problem \eqref{splines exp. on t and r, opt.problem} is 
 \begin{equation*}
    argmin_{\theta'\in \Theta}\frac{1}{3M}\frac{1}{n_{t}}\frac{1}{n_{p}}\sum_{n=1}^{n_p}\sum_{i=1}^{n_{t}}\sum_{I=1}^{M}||F_{I}(q_{i,n})- \sum_{s=1}^{N_{d}\times N_{b}}\theta'_{s}\cdot J_{I;s}(\COP q_{i,n})||^{2}\PERIOD    
\end{equation*}
A matrix form can represent the last expression as
\begin{equation}
argmin_{\theta\in\Theta}\frac{1}{3Mn_{t}n_{p}}\left(\mathbf{F}-\mathbf{J}\theta'\right)^{tr}\left(\mathbf{F}-\mathbf{J}\theta'\right)\COMMA\label{otimization pr with matrix new}
\end{equation}
where $\mathbf{F}$ is a $3Mn_{t}n_{p}\times 1$ vector, $\mathbf{J}$ is a $3Mn_{t}n_{p}\times N_{d}N_{b}$ matrix and $\theta'$ is a $N_{d}N_{b}\times 1$ vector. Finding the solution of the optimization problem \eqref{otimization pr with matrix new} is equivalent to finding that parameter set $\theta'$ that solves

\begin{equation}
\mathbf{J}\theta'=\mathbf{F}\COMMA\label{normal eq}
\end{equation} 
where

 $$  \mathbf{J}=
\begin{bmatrix}
  J_{1x;1}(\COP q_{1,1})
  
   &   J_{1x;2}(\COP q_{1,1})
 
   &\dots 
   
   &   J_{1x;N_{d}N_{b}}(\COP q_{1,1})
 \\
    J_{1y;1}(\COP q_{1,1})
  & 
   J_{1y;2}(\COP q_{1,1})
    
   & \dots  
    &
     J_{1y;N_{d}N_{b}}(\COP q_{1,1})
 \\
   \vdots  & \vdots 
   &\dots  &\vdots \\
      J_{Mz;1}(\COP q_{n_{t},n_{p}})
  &  J_{Mz;2}(\COP q_{n_{t},n_{p}})
 &\dots  &   J_{Mz;N_{d}N_{b}}(\COP q_{n_{t},n_{p}})
\end{bmatrix}
,$$
 with the elements, $J_{i;j} $ of the matrix $ \mathbf{J}$  are the values of the representation at the observations, for example, ${J_{1x;1}(\COP q_{1,1})=\chi_{1}(t)\cdot\sum_{J>1}^{M}\psi_{1}(r_{1J})\frac{\partial r_{1J}}{\partial x}, \dots,}\\  J_{1x;NdNb}(\COP q_{1,1})=\chi_{Nb}(t)\cdot\sum_{J>1}^{M}\psi_{Nd}(r_{1J})\frac{\partial r_{1J}}{\partial x}  $.

 \subsection{Estimation of the friction kernel}\label{representation of the friction kernel subsection}
We estimate the dynamic friction kernel from the underlying AA MD simulations using correlation functions of the CG particles \cite{VothDyn2006,Lei14183}. 
A more realistic form of dynamics \eqref{eq:GLE} is 
\begin{equation}
   \delta \bar{F}(\bx_{t}) =-\int_{0}^{t}\bar{\zeta}(t-\tau)\bar{v}_{\tau}d\tau+F_{R}(t)\COMMA\label{GLE2}
\end{equation}
where $\delta\bar{F}(\bx)=  F_{tot}(\bx)-\FPMF(\bx) $ is the instantaneous difference between the mean CG forces $\FPMF(\bx)$ and the exact AA MD forces $F_{tot}(\bx_{t})$, and $F_{R}(t)$ denotes the random force, and $\bar{v} = \masscg^{-1}\bp$. Multiplying both sides of \eqref{GLE2} by the initial CG velocity vector, $\bar{v}(0)$, and taking the canonical average, we obtain the following relation,

\begin{equation}
   C_{\delta \bar{F}\bar{v}}(t)=-\int_{0}^{t}\bar{\zeta}(t-\tau)C(\tau)d\tau\COMMA\label{volterra2}
\end{equation}
where  $C(t)=\frac{1}{3}\langle \bar{v}(t)\cdot\bar{v}(0)\rangle$, and $C_{\delta \bar{F}\bar{v}}(t)=\frac{1}{3}\langle \delta \bar{F}(t)\cdot\bar{v}(0)\rangle$ are the velocity auto-correlation function and the force-velocity correlation function of the CG particles, respectively which can be directly evaluated from the AA MD simulations, while $\langle\cdot\rangle$, is defined on the probability space $(\Omega, \Ff, \P)$. Considering that the initial velocity and random force, $F_{R}$, are uncorrelated, the corresponding correlation function $\langle F_{R}(t)\cdot\bar{v}(0)\rangle$ vanishes. Taking the Laplace transform on both sides of (\ref{volterra2}), we obtain,
\begin{equation}
    \Lapl{C}_{\delta \bar{F}\bar{v}}(s)=-\Lapl{\bar{\zeta}}(s)\Lapl{C}(s)\COMMA\label{EQ2}
\end{equation}
where $\Lapl{C}(s),\:\Lapl{C}_{\delta \bar{F}\bar{v}}(s)$ and $\Lapl{\bar{\zeta}}(s)$ denote the Laplace transform of $C(t)$, $C_{\delta \bar{F}\bar{v}}(t)$ and $\bar{\zeta}(t)$ respectively. We use the hat $"\ \hat{\ }\ "$ notation to denote the Laplace transform in this subsection. 

Solving \eqref{EQ2} for $\Lapl{\bar{\zeta}}(s)$, we derive an equation that relates the dynamic friction kernel with the correlation functions in the Laplace domain,

\begin{equation}
  \Lapl{\bar{\zeta}}(s)=-\frac{\Lapl{C}_{\delta \bar{F}\bar{v}}(s)}{\Lapl{C}(s)}\PERIOD\label{relation friction with CFs 2}
\end{equation}
In the limit, $s=0$, the Laplace transform of the friction kernel is expressed as 
\begin{equation}
  \Lapl{\bar{\zeta}}(0)=-\frac{\Lapl{C}_{\delta \bar{F}\bar{v}}(0)}{\Lapl{C}(0)}=-\frac{\frac{1}{3}\int_{0}^{\infty}\langle \delta \bar{F}(t)\cdot\bar{v}(0)\rangle dt}{\frac{1}{3}\int_{0}^{\infty}\langle \bar{v}(t)\cdot\bar{v}(0)\rangle dt}\PERIOD\label{friction coeff formula}
\end{equation}
The two quantities $\int_{0}^{\infty}\langle \delta \bar{F}(t)\cdot\bar{v}(0)\rangle dt $ and $\int_{0}^{\infty}\langle \bar{v}(t)\cdot\bar{v}(0)\rangle dt$ are Green-Kubo-type of relations. Green-Kubo is an important class of relations in which a macroscopic dynamical property is written with respect to microscopic time correlation functions \cite{Evans-Morris}.

It is convenient to approximate the random force as a Markovian force (as is the case in dynamics \eqref{eq:GLE}) such that 

 \begin{equation}
        \langle F_{R}(t)\cdot F_{R}(0)\rangle\approx 2\beta^{-1}\bar{\zeta}_{0}\delta(t),
    \end{equation}
where $\bar{\zeta}_{0}$ is a constant and $\delta(t)$ is the dirac function. Choosing $\bar{\zeta}_{0}=\tilde{\bar{\zeta}}(0)$, the friction coefficient is expressed as 

\begin{equation}
 \bar{\zeta}_{0}=-\frac{\frac{1}{3}\int_{0}^{\infty}\langle \delta \bar{F}(t)\cdot\bar{v}(0)\rangle dt}{\frac{1}{3}\int_{0}^{\infty}\langle \bar{v}(t)\cdot\bar{v}(0)\rangle dt}\PERIOD\label{zeta0}
\end{equation}

In this work, we estimate the friction coefficient using two different atomistic data sets. In the first case, the AA data set,  $\left\{ v_{t},f_{t} \right\}_{a\leq t\leq b}$  of atomistic velocities and forces, corresponds to the equilibrium regime, $[a,b]$. Thus, the empirical estimator of \eqref{zeta0} is

\begin{equation}\label{eq: equilibrium friction formula}
 \fe= \bar{\zeta}_{0}=-\frac{\frac{1}{3}\int_{a}^{b}\langle \delta \bar{F}(t)\cdot\bar{v}(0)\rangle dt}{\frac{1}{3}\int_{a}^{b}\langle \bar{v}(t)\cdot\bar{v}(0)\rangle dt}\COMMA
\end{equation}
where $\delta\bar{F}(t)=  F_{tot}(\bx_{t})-\bar{F}(\bx;\phi^{\ast}) $ is the instantaneous difference between the effective equilibrium CG forces $\bar{F}(\bx;\phi^{\ast})$ as derived by solving the optimization problem \eqref{discrete PSFM} and the exact AA MD forces $F_{tot}(\bx_{t})$.

 In the second case, we estimate the friction coefficient using transient AA information. Our data set, $\left\{ v_{t},f_{t} \right\}_{  t\in \tr}$, corresponds to the system's transient regime, $\tr$, where the system's initial state is out of equilibrium. The empirical estimator is
\begin{equation}\label{eq: transient friction formula}
 \ft= \bar{\zeta}_{0}=-\frac{\frac{1}{3}\int_{0}^{t_{f}}\langle \delta \bar{F}(t)\cdot\bar{v}(0)\rangle dt}{\frac{1}{3}\int_{0}^{t_{f}}\langle \bar{v}(t)\cdot\bar{v}(0)\rangle dt}\COMMA
\end{equation}
where $\delta\bar{F}(t)=  F_{tot}(\bx_{t})-\bar{F}(\bx,t;\theta^{\ast}) $ is the instantaneous difference between the effective time-depended CG forces $\bar{F}(\bx,t;\theta^{\ast})$ as derived by solving the optimization problem \eqref{splines exp. on t and r, opt.problem} and the exact AA MD forces $F_{tot}(\bx_{t})$. 
\section{Benchmark example: one particle described via the Generalized Langevin Equation}\label{benchmark}
We first consider a relatively simple model concerning the dynamics of one moving particle in a box that is described by the process $\{\bx_{t},\bp_{t}\}_{t\geq0}$, $(\bx, \bp)\in \mathbb{R}^{2}$ satisfying the GLE

\begin{equation}\label{eq:GLE_2}
    \begin{cases}
    dQ_t = P_t dt , \\
   dP_{t} = -\alpha Q_tdt - [\eta^2\,\int_0^t   e^{\frac{t-s}{\tau} }P_s ds] + F_{R}({t}),\quad t>0,\\   (Q_{0},P_{0})\sim \nu_0(d\bx,d\bp),
    \end{cases}
\end{equation}
where the particle is further subjected to a harmonic oscillator potential $U(\bx)=\frac{1}{2}\alpha\bx^{2}$, and
$\alpha$ reflects the strength of the oscillator. The quantity $\bar{\zeta}(t)=\eta^2e^{\frac{t}{\tau} }$ is the friction kernel, while $\eta$ is a constant controlling the amplitude of the friction kernel.

The auto-correlation function of the random force $F_{R}(t)$ is given by 
 \begin{equation} \label{eq:correlation_function}
  C_{F_{R}}(t-s)=\beta^{-1}\, \eta^2 e^{\frac{t-s}{\tau}}, \quad t>s.   
 \end{equation}
 where $\tau$ is a constant which controls the rate of decay of \eqref{eq:correlation_function} \cite{pavliotis2014stochastic}.

We solve numerically  the equivalent extended dynamics  of \eqref{eq:GLE_2}, \cite{pavliotis2014stochastic},  to generate   approximate sample paths $\{Q_{t_i},P_{t_i}\},i=1,...,n_t$ of $Q_{t}, P_{t}$ using the Euler-Maruyama numerical scheme \cite{kloeden2013numerical}. 
From the obtained data, we infer the parameters of the model Langevin equation with a time-dependent   force  $  F_{CG}(\tilde{Q},t;\theta)$ 
\begin{equation}
\begin{cases} 
     d\tilde{Q}_t = \tilde{P}_t dt \\
    d\tilde{P}_t = F_{CG}(\tilde{Q}_t,t;\theta) -\bar{\zeta}_{0}   \tilde{P}_t dt+ \bar{\sigma}_{0} d\bar{W}_t, \quad t>0, \\
   ( \tilde{Q}_{0} \tilde{P}_{0}) \sim \bar\nu_0
   (d\tilde{\bx},d\tilde{\bp})\PERIOD
 \end{cases} \label{eq:time_d_potential}
 \end{equation}
Let us denote  the total force 
 \begin{equation*}
  F_{tot}(Q_t,\bp_{t},t)= -\alpha Q_t - [\eta^2\,\int_0^t  e^{\frac{t-s}{\tau} }P_s ds] + F_{R}({t})\PERIOD
 \end{equation*}
We choose  the parameterized representation of the CG forces with respect to $t$ and $\tilde{Q}$  
 \begin{equation*} 
  F_{CG}( {Q},t;\theta)=-D(t;\theta^{1}) \,B( {Q} ;\theta^{2})\COMMA
 \end{equation*}
  where $\theta=\{ \theta^{1}, \theta^{2}\}$ is the  parameter set. 

To find the optimal parameter set $\theta^{\ast}$, we solve the PSFM optimization problem given by
  \begin{equation}
    \theta^{\ast}_{\tr}= argmin_{\theta\in \Theta}  \frac{1}{n_{t}}\frac{1}{n_{p}}\sum_{n=1}^{n_p}\sum_{i=1}^{n_{t}} ||F_{tot}(Q_{t_{i}},P_{t_{i}},t_{i})-F_{CG}(Q_{t_{i}},t_{i};\theta)||^2 \COMMA\label{eq:FM_tdp_1}
 \end{equation}
 where $n_p$ is the total number of sample paths, and $n_t$ is the number of data points in each path. $||\cdot||$ denotes the Euclidean norm in $\R$. The data are sampled in $\tr$ from the {\it{initial for PSFM}} data set, see the tables \ref{datasettau05} and \ref{datasettau01}. 
 
We invoke the FDT to estimate $\bar{\zeta}_{0} $ from $\bar{\sigma}_{0}$ using the relationship $\bar{\sigma}^2_{0}= {2\, \beta^{-1}\bar{\zeta}_{0} }$.
Then, for the $\bar{\sigma}_{0}^2$ we use as  a consistent estimator the quadratic variation of the process \cite{iacus2008simulation,pavliotis2008}
 \begin{equation}
      \bar{\sigma}_{0}^2= \frac{1}{n_t^{\sigma}\Delta}\sum_{i=1}^{n_{t}^{\sigma}}(P_{t_{i}}-P_{t_{i-1}})^2 \label{eq:sample_variance}, \quad \Delta = t_i-t_{i-1},
 \end{equation}
 where the $n^{\sigma}_t$ data points of the {\it{diffusion coefficient data set}} 
 are sampled by an {\it initial} data set
 obtained by the Euler-Maruyama scheme from one path that is long enough to include data of both the transient and stationary time regimes. 
 For the estimation of $\bar{\sigma}_{0}$, $\bar{\zeta}_{0} $ of the time-dependent Langevin equation, we consider a sparser grid of the long path. This coarse approach in the time-domain \cite{bernardes2017coarse} is connected to the non-Markovian character of the actual system. Sampling the system's evolution in time steps that are longer than the system's correlation time, as defined below, makes the Markovian property a more plausible assumption and allows us to obtain an effective diffusion coefficient for the proposed Markovian model.
  
\begin{table}[htbp]
\resizebox{\textwidth}{!}{
\begin{tabular}{|c|c|c|c|c|}
\hline
Data sets        & Paths $n_p$ & Observation Frequency & Observations per path  & Observation  time intervals $\tr$ \\ \hline
Initial for PSFM         & 2000        & 0.005                 & 1000                        & $[0,5]$                     \\ \hline
PSFM             & 2000        & 0.05                  & 100                         & $[0,5]$                     \\ \hline
Initial for diff. coef. & 1           & 0.15                  & 1000                        & $[0,150]$                   \\ \hline
Diff. coef.     & 1           & 0.83                  & 180                         & $[0,150]$                   \\ \hline
\end{tabular}}
\\
\caption{Details of the data sets used for the solution of the PSFM problem for the correlation time $\tau=0.5$.}
\label{datasettau05}
\end{table}
   
\begin{table}[htbp]
\resizebox{\textwidth}{!}{
\begin{tabular}{|c|c|c|c|c|}
\hline
Data sets        & Paths $n_p$ & Observation Frequency & Observations per path  & Observation  time intervals $\tr$\\ \hline
Initial for PSFM         & 2000        & 0.012                 & 1000                        & $[0,12]$                     \\ \hline
PSFM             & 2000        & 0.12                  & 100                         & $[0,12]$                     \\ \hline
Initial for diff. coef. & 1           & 0.15                  & 1000                        & $[0,150]$                   \\ \hline
Diff. coef.     & 1           & 0.1875                  & 800                         & $[0,150]$                   \\ \hline
\end{tabular}}
\caption{Details of the data sets used for the solution of the PSFM problem for the correlation time $\tau=0.1$.}
\label{datasettau01}
\end{table}
To check the validity of our approach with respect to the non-Markovian character of the GLE system, we present results for two different values of the parameter $\tau$ of the GLE. The parameter $\tau$ controls the correlation time $t_F^{corr}=\int_0^\infty C_{F^{rn}}(\tau) d\tau/ C_{F^{rn}}(0)=\tau$ of the colored noise $F_{R}(t)$. We consider a more correlated force for  $\tau=0.5$ and a less correlated force  for $\tau=0.1$. The other parameters of the GLE  are $\alpha=0.1,\eta=0.1, \beta=1$, $Q_{0}=0, P_{0}=0.1$. The characteristics  of the data sets used for the solution of the PSFM problem are listed in table \ref{datasettau05} for $\tau=0.5$ and in table \ref{datasettau01} for $\tau=0.1$.  
The {\it initial for PSFM} data set and the {\it PFSM} data set  have  the observation time interval $t\in [0,5]$ for $\tau=0.5$  and $t \in [0,12]$ for $\tau=0.1$. The {\it initial for PSFM} data set
consists of $n_p=2000$ paths and $1000$ data points.  From every path of the {\it initial for PSFM} data set, we pick  $n_t=100$ data points to obtain the {\it PSFM} data set for both correlation times. 

For the representation of $D(t;\theta^{1})$ and $B( {Q} ;\theta^{2})$, we use a non-parametric representation with cubic splines with $N_b$ and $N_d$ basis functions, respectively. For $\tau=0.5$ we assume that  $N_b=N_d=10$ whereas for  $\tau=0.1$ we assume that  $N_b=N_d=15$.    

We estimate the diffusion coefficient using sub-sampling, that is from the {\it diffusion coefficient} data set that is a subset  of the {\it{initial for diffusion coefficient}} data set. The data sets used for the diffusion coefficient have a longer observation interval with  $t\in [0,150]$ and consist of one path, as listed, along with their other characteristics in table~\ref{datasettau05} for $\tau=0.5$  and in  table~\ref{datasettau01} for $\tau=0.1$.  For $\tau=0.5$, according to \eqref{eq:sample_variance}, we  estimated  $\bar{\sigma}_{0}=0.097$   using $n_{t}^{\sigma}=180$ data points that correspond to a time step (observation frequency) of $0.83$. For $\tau=0.1$, we estimated $\bar{\sigma}_{0}=0.044$ using $n_{t}^{\sigma}=800$ data points that correspond to a time step (observation frequency) of $0.1875$. 

\begin{figure}[htbp]
\centering
\subcaptionbox{}[.49\textwidth]{\includegraphics[width=0.49\textwidth]{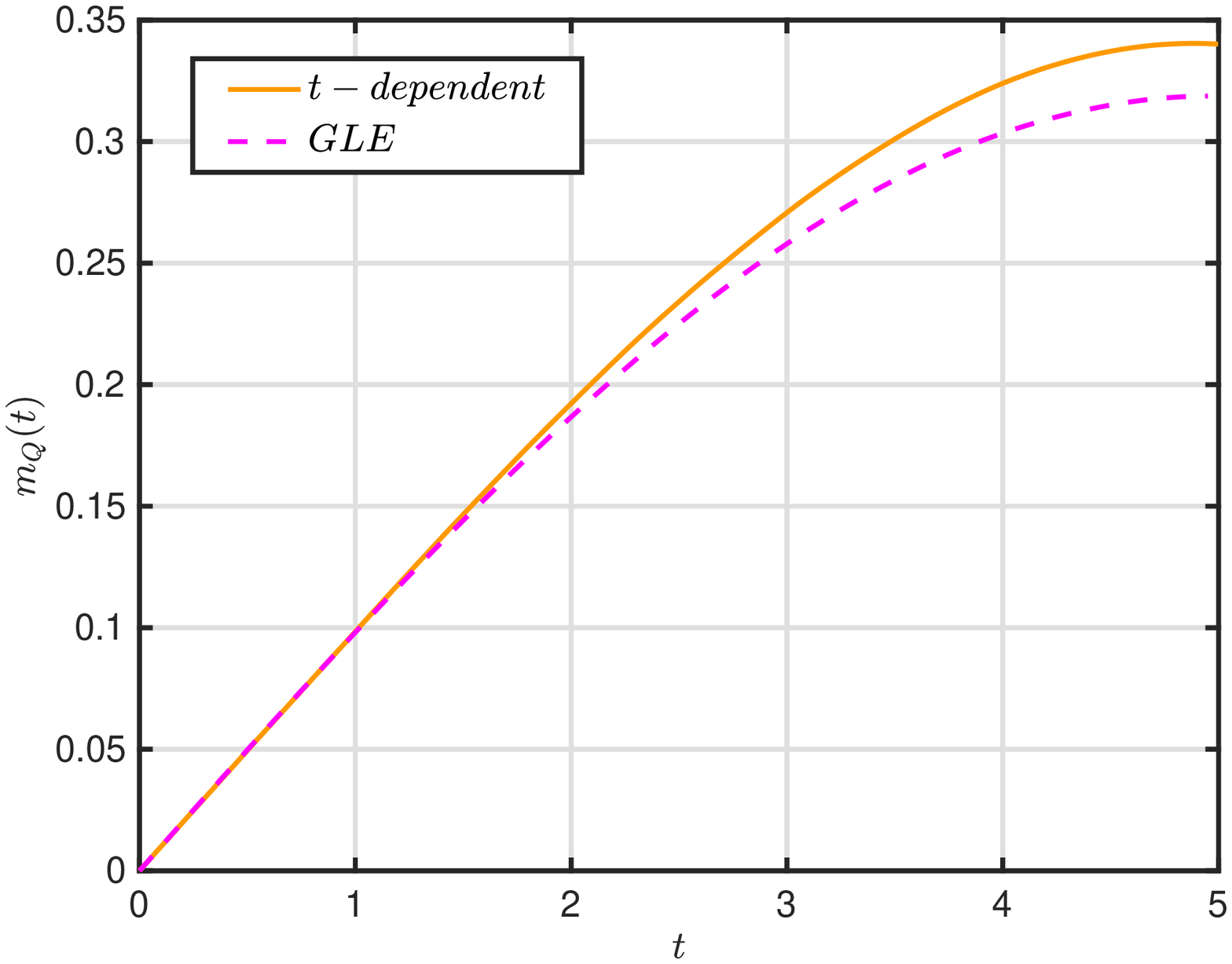}}
\subcaptionbox{}[.49\textwidth]{\includegraphics[width=0.49\textwidth]{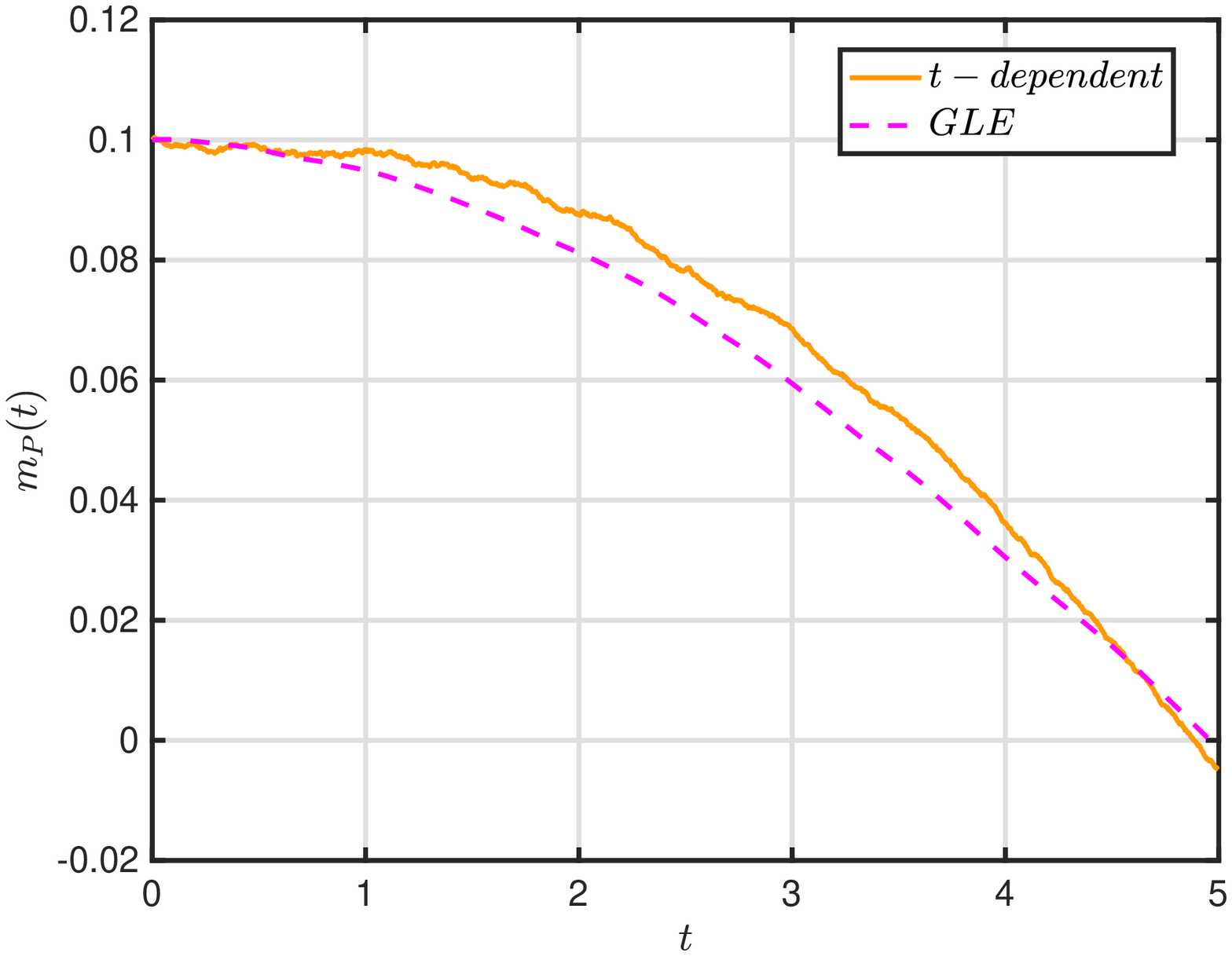}}\\
\subcaptionbox{}[.49\textwidth]{\includegraphics[width=0.49\textwidth]{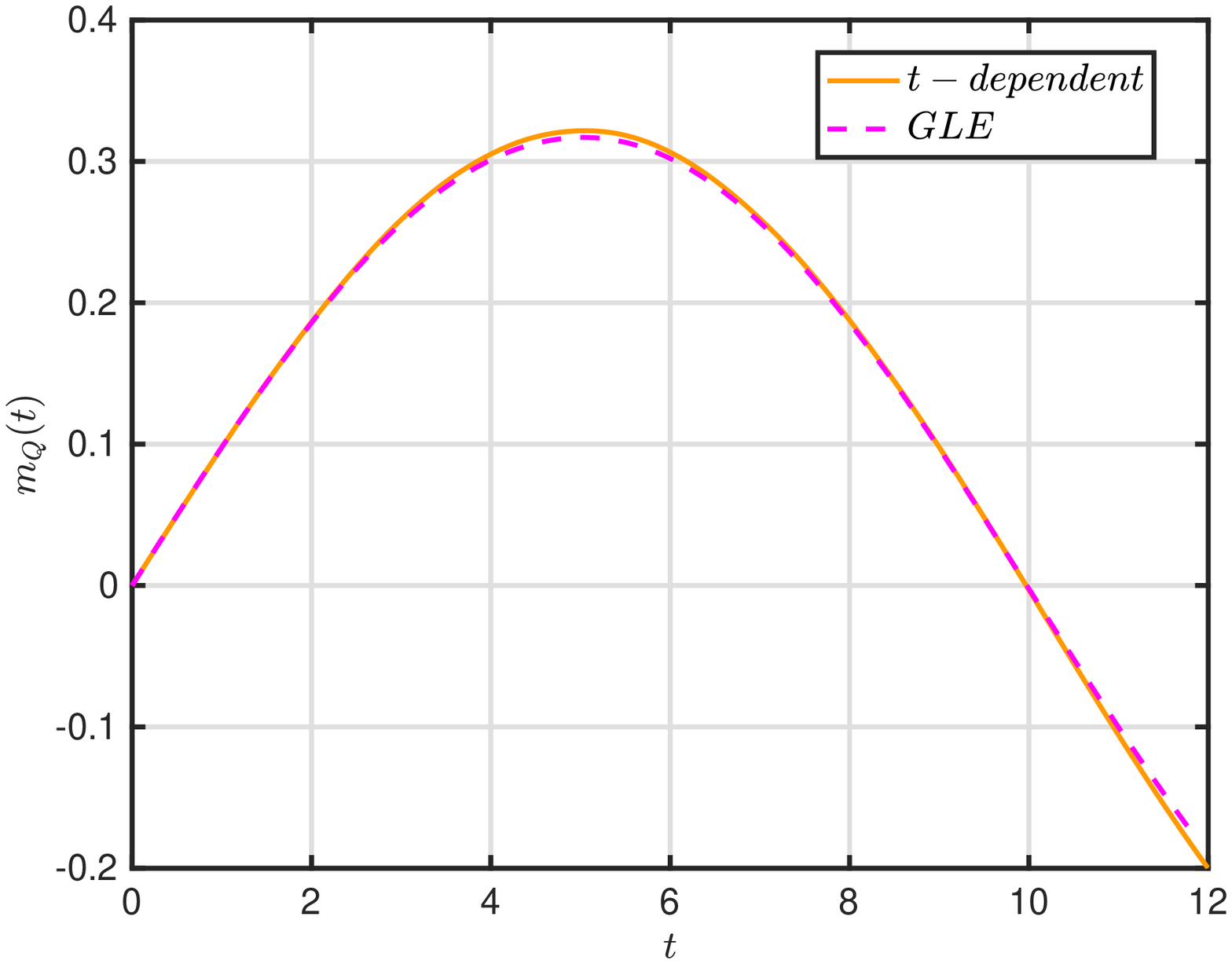}}
\subcaptionbox{}[.49\textwidth]{\includegraphics[width=0.49\textwidth]{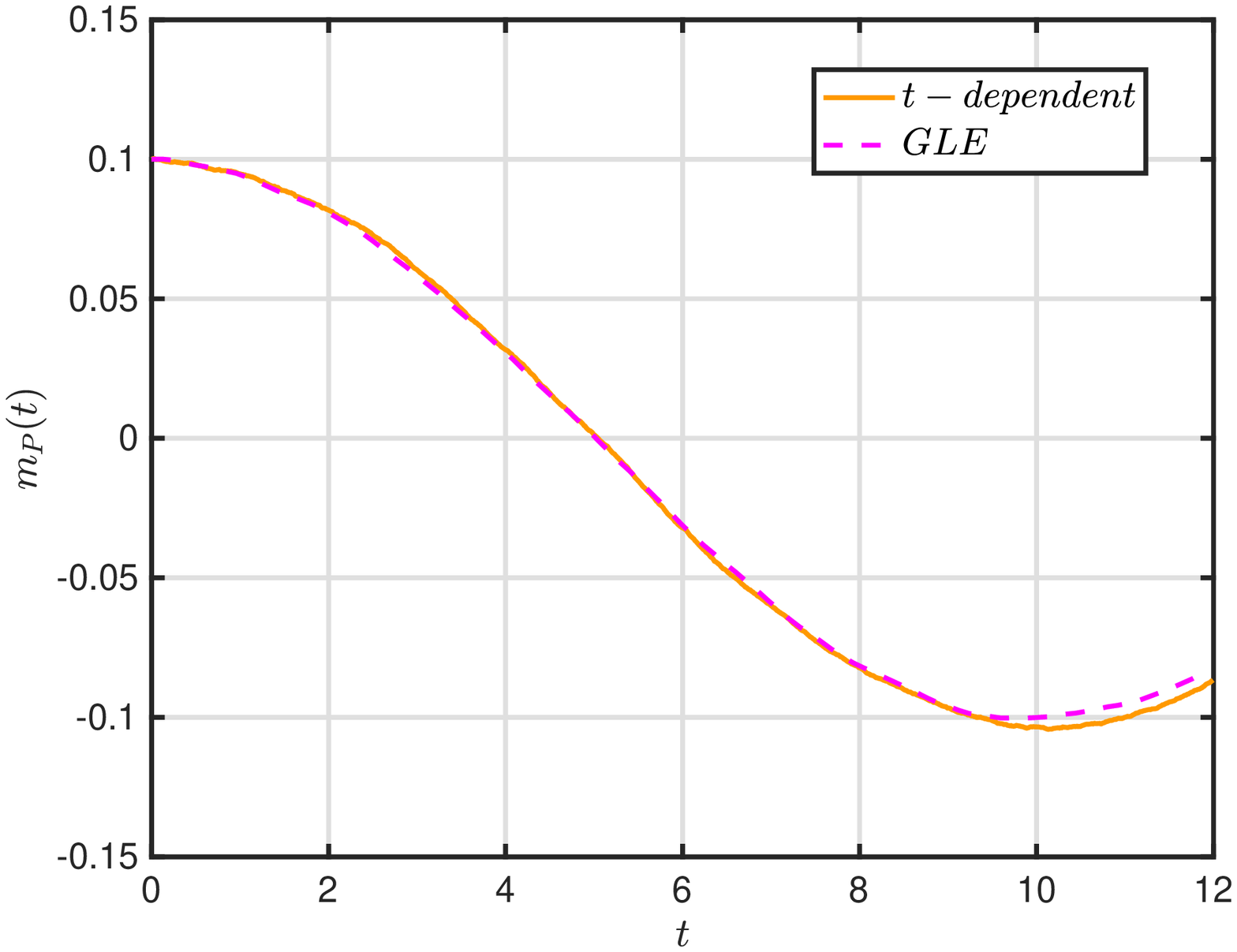}}
\caption{ 
The sample mean values of the GLE and of the Langevin equation, for a single particle model, with a time-dependent potential. 
For correlation times $\tau=0.5$ (a,b) and $\tau=0.1$ (c,d). A longer interval of the transient state of the  GLE system dynamics is well approximated by the Langevin equation with time-dependent potential for smaller correlation time. 
 }
\label{fig:GLE_mean}
\end{figure}

\begin{figure}[htbp]
\centering
\subcaptionbox{}[.49\textwidth]{\includegraphics[width=0.49\textwidth]{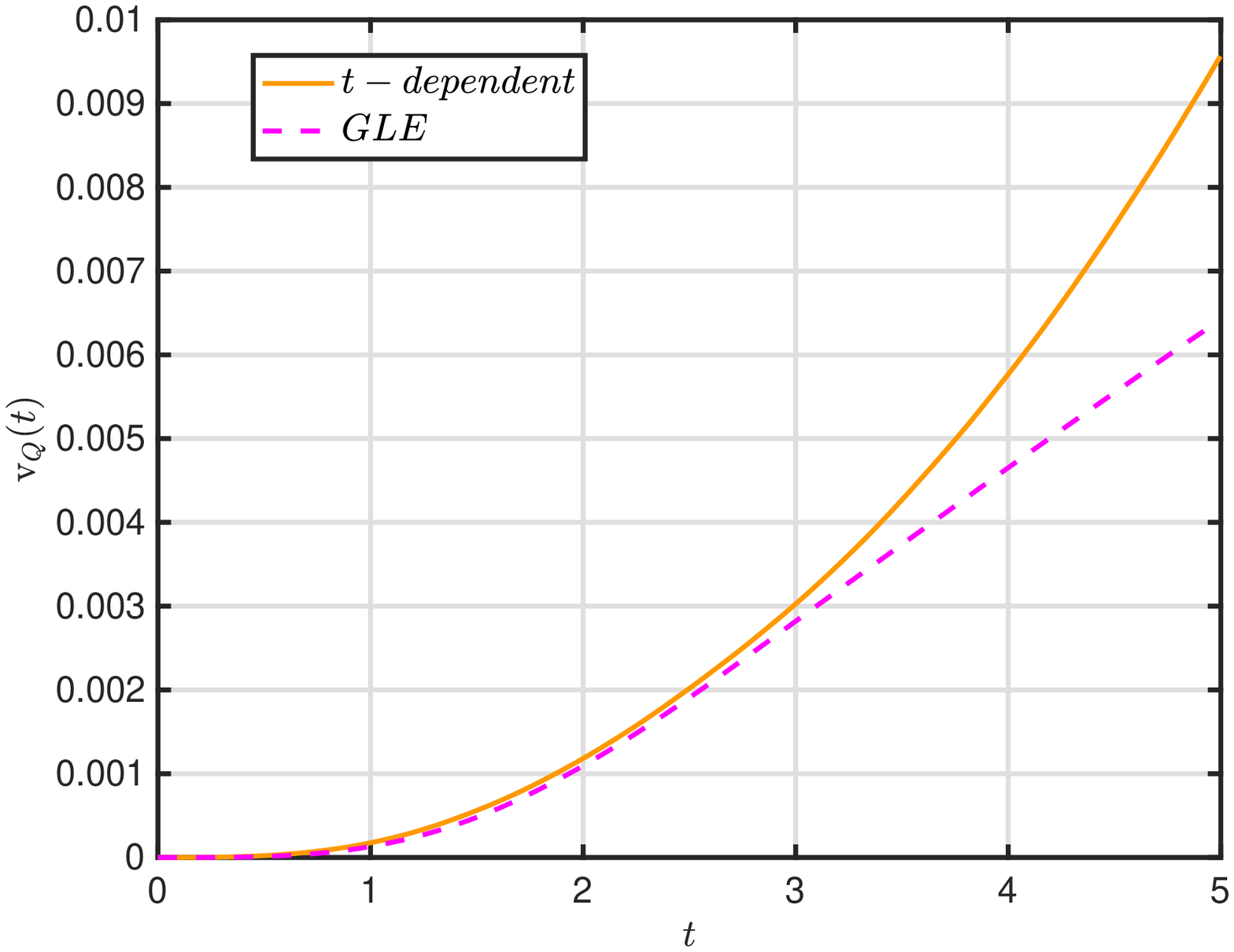}}
\subcaptionbox{}[.49\textwidth]{\includegraphics[width=0.49\textwidth]{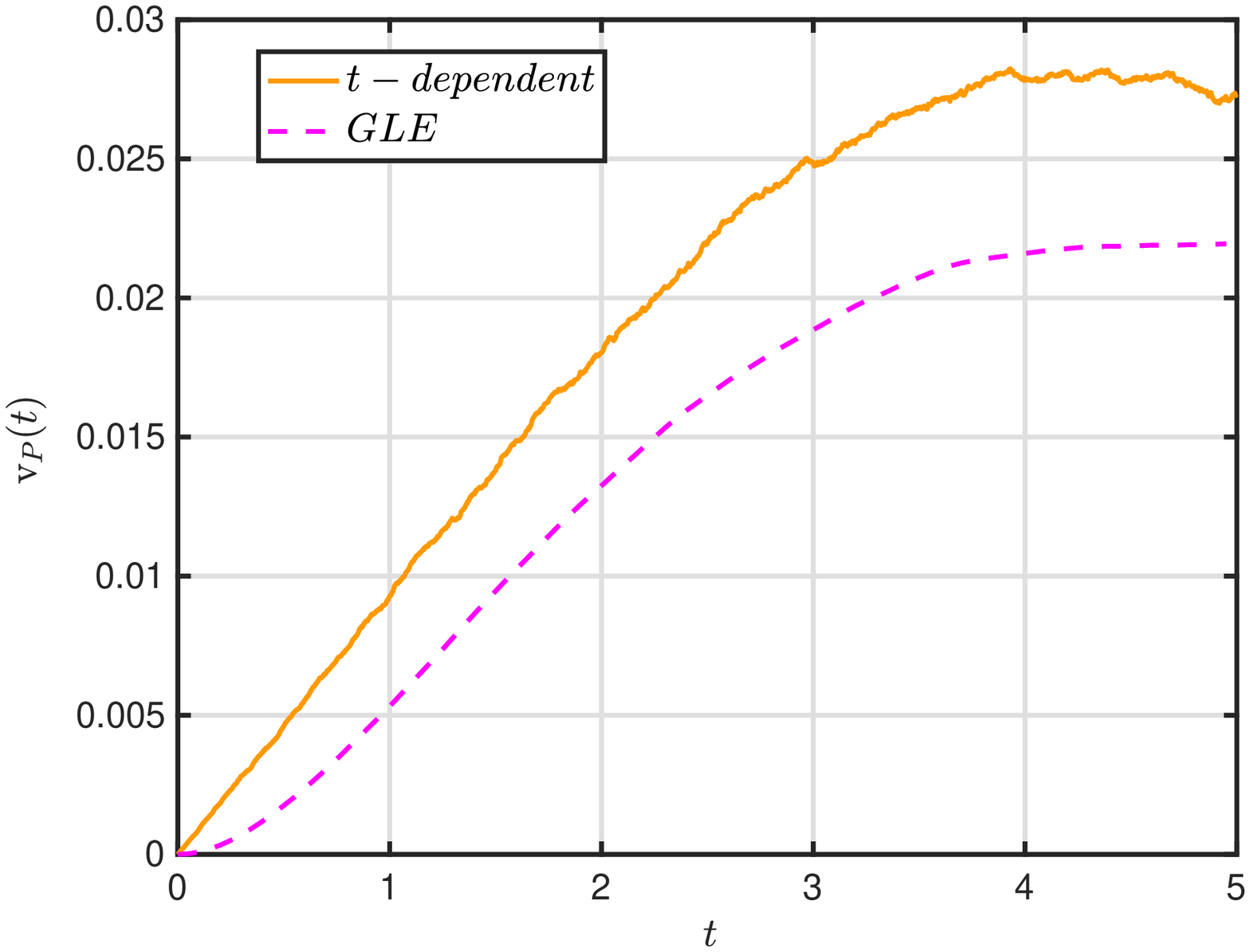}}\\
\subcaptionbox{}[.49\textwidth]{\includegraphics[width=0.49\textwidth]{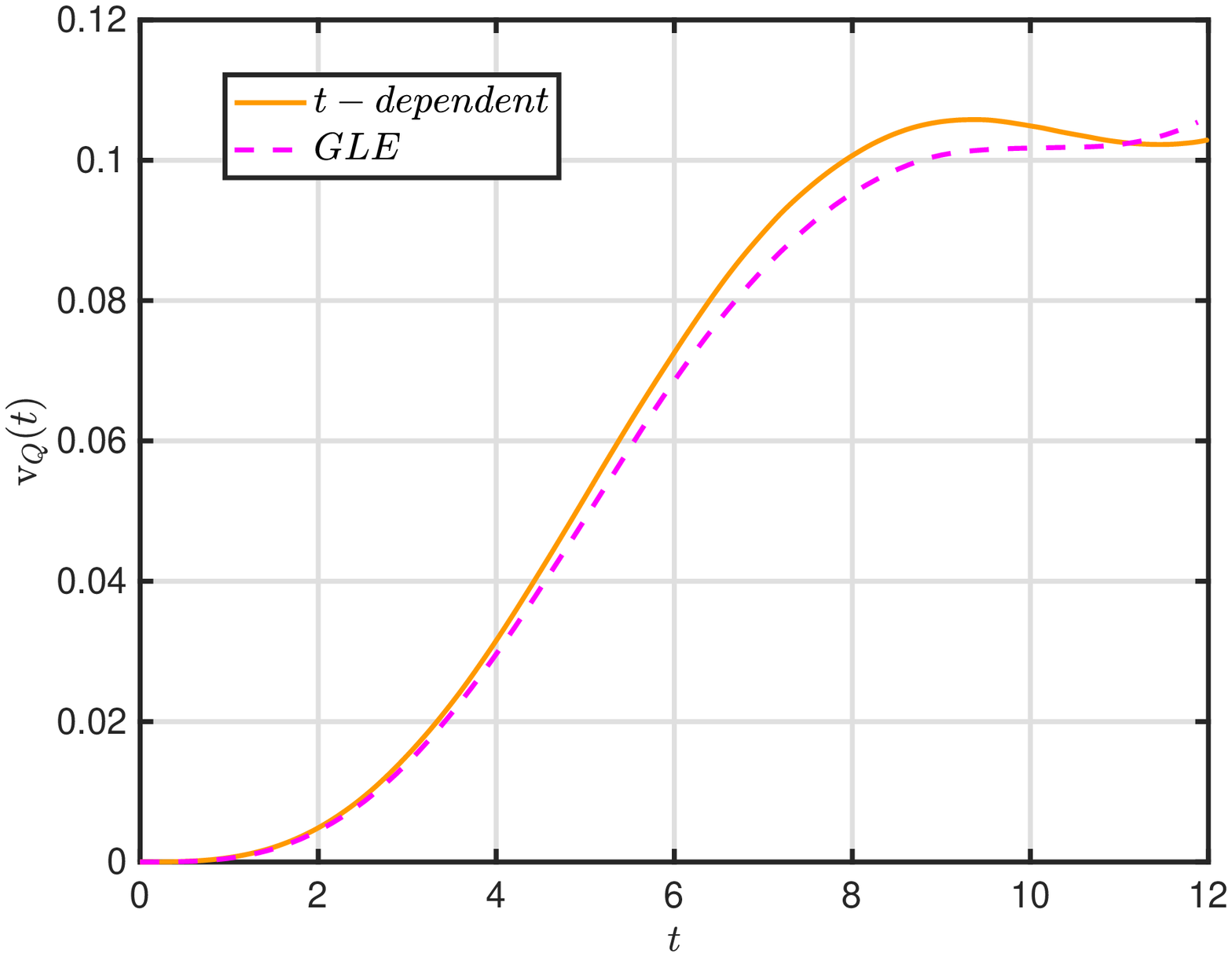}}
\subcaptionbox{}[.49\textwidth]{\includegraphics[width=0.49\textwidth]{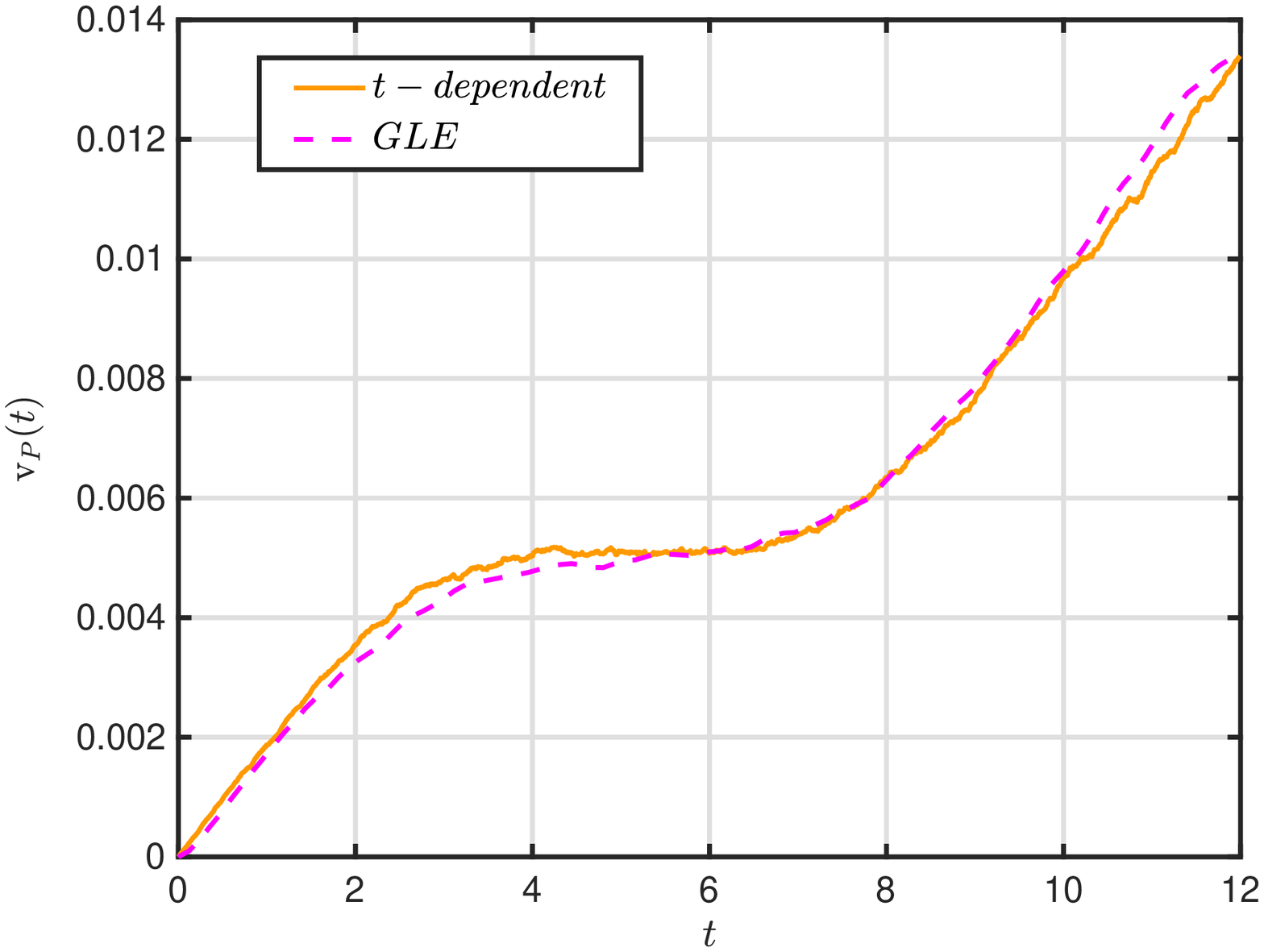}}
\caption{Comparison of the sample variances of the GLE and of the Langevin equation, for a single particle model, with a time-dependent potential.} \label{fig:GLE_variance}
\end{figure}

 In Fig.\ref{fig:GLE_mean} and  Fig.~\ref{fig:GLE_variance}, we compare the empirical  mean value and variance of  $Q_{t},P_{t}$ respectively of the GLE  (dashed lines) with the corresponding of $\tilde{Q}_t, \tilde{P}_t$  of the  Langevin equation with time-dependent potential (solid lines) whose parameters were estimated using the scheme described above. 
 In the upper panel, results correspond to the more correlated case for $\tau=0.5$, whereas the lower panel results correspond to the less correlated case for $\tau=0.1$. We notice that the time-dependent Langevin equation can capture the transient short-time regime for both cases. However, the model is more efficient for the less correlated case for $\tau=0.1$ when the actual system is closer to a Markovian system. In the less correlated case, the mean values of $Q_{t},P_{t} $ are well approximated for a longer time interval.

It is important to point out that the results for the variances are sensitive to the correct estimation of $\bar{\sigma}_{0}, \bar{\zeta}_{0} $ from the data. For both cases presented herewith, we chose a longer path to include data from both the transient-time and the steady-state regime. Moreover, we used a sparser grid for the more correlated case than the less correlated case. This coarse approach in the time domain prevents the under-estimation of  $\bar{\sigma}_{0}$, $\bar{\zeta}_{0}$, which would result in a significant under-estimation of the variance of $Q_t, P_t$. 
We expect that the system dynamics can be better approximated for a longer interval of the transient-time regime when time-dependent $\bar{\sigma}_{0},\bar{\zeta}_{0}$ are included in the time-dependent Langevin Equation \eqref{eq:time_d_potential}. 

Summarizing, with this benchmark example  we have verified that the proposed model can approximate well the short transient time regime of the GLE describing the movement of one particle in a box. The   correlation time $\tau$ is the parameter that controls the range of validity of the model.  

\section{Transient time dynamics for a liquid water system}\label{results}

\begin{figure}[htbp]
    \centering
    \includegraphics[width=0.5\textwidth]{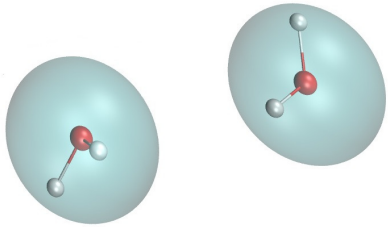}
    \caption{Atomistic and CG representations of two water molecules.}
    \label{cg mapping}
\end{figure}

Our following example involves a  realistic molecular system, that of a liquid water droplet consisting of $N$ molecules. 
We examine a system of $M$ water molecules where a CG particle is the center of mass of a group of atoms that correspond to a single water molecule, see Fig. \ref{cg mapping}. 
Thus, there are only non-bonded pair interactions between the CG particles, and the effective CG interaction potential is sought in the form,
\begin{equation*}
  \bar{U}(\bx)=\sum_{I=1}^{M}\sum_{J>I}^{M} \bar{u}(r_{IJ})\COMMA
\end{equation*} 
where $\bar{u}(r)$ is the potential that corresponds to the pair interaction between CG particles $I$ and $J$, at distance $r_{IJ}$, such that $r_{IJ}=|\bx_{I}-\bx_{J}|,\:I\neq J,\:I,J=1,..,M.$ 

We impose the initial, far from equilibrium, state of the water system to be an artificial FCC (Face Centered Cubic) configuration, of a very low density, so that we are able to examine the transition towards the chosen thermodynamic equilibrium of the system, i.e. a water nanodroplet, see Fig. \ref{fcc initial snapshot}. For this system, we estimate the transient dynamics described by the proposed model in \eqref{eq:Langevin_time_d}. We also perform equilibrium MD simulations to compare the proposed CG model with equilibrium dynamics models. We further compare the proposed CG dynamics characterized by the time-dependent CG potential with the Langevin dynamics given by \eqref{eq:Langevin equation}. To obtain the complete dynamics models we estimate the friction coefficient for the proposed dynamics in two ways. In the first, we use data corresponding to the short-time transient regime, while in the second, we use data from equilibrium, see section~\ref{representation of the friction kernel subsection}. 

The results presented in this section (a) confirm that the widely used Markovian models do not apply in the transient short-time dynamics and (b) establish the validity and consistency of the proposed dynamic models. 

\subsection{Atomistic simulations and configuration data sets}
\begin{figure}[htbp]
 \centering
\subcaptionbox{}[.32\textwidth]{\includegraphics[width=0.32\textwidth]{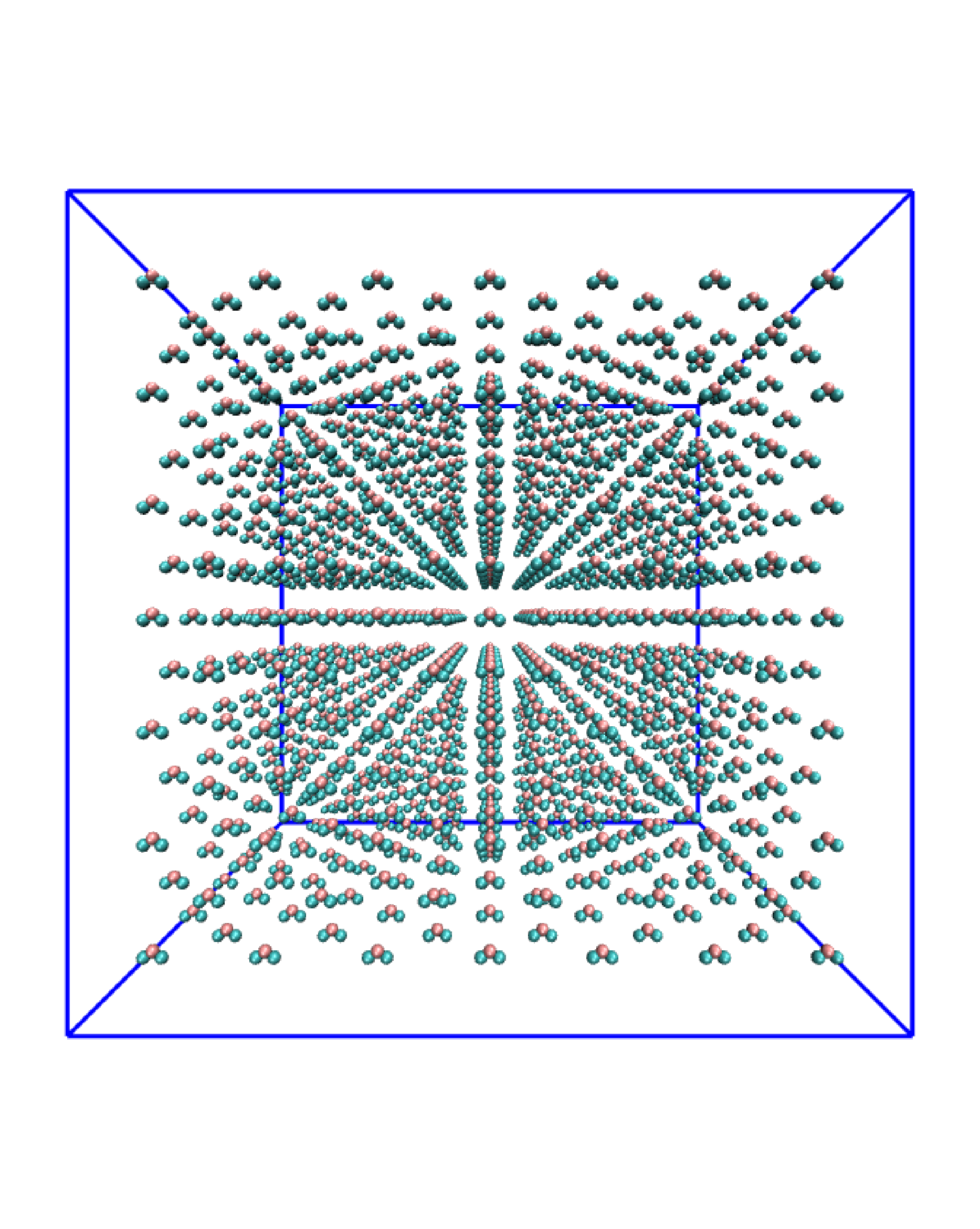}}
\subcaptionbox{}[.32\textwidth]{\includegraphics[width=0.32\textwidth]{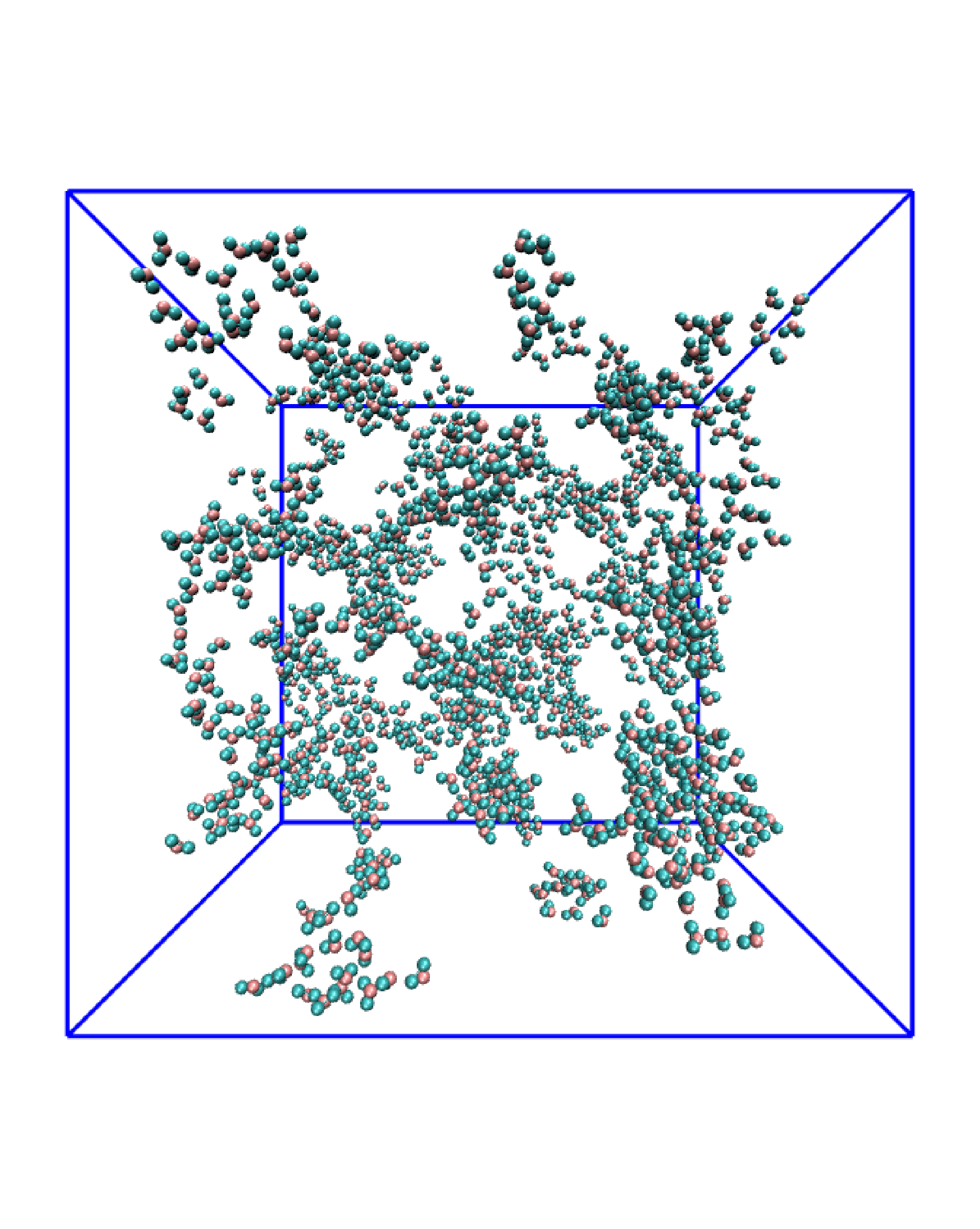}}
\subcaptionbox{}[.32\textwidth]{\includegraphics[width=0.32\textwidth]{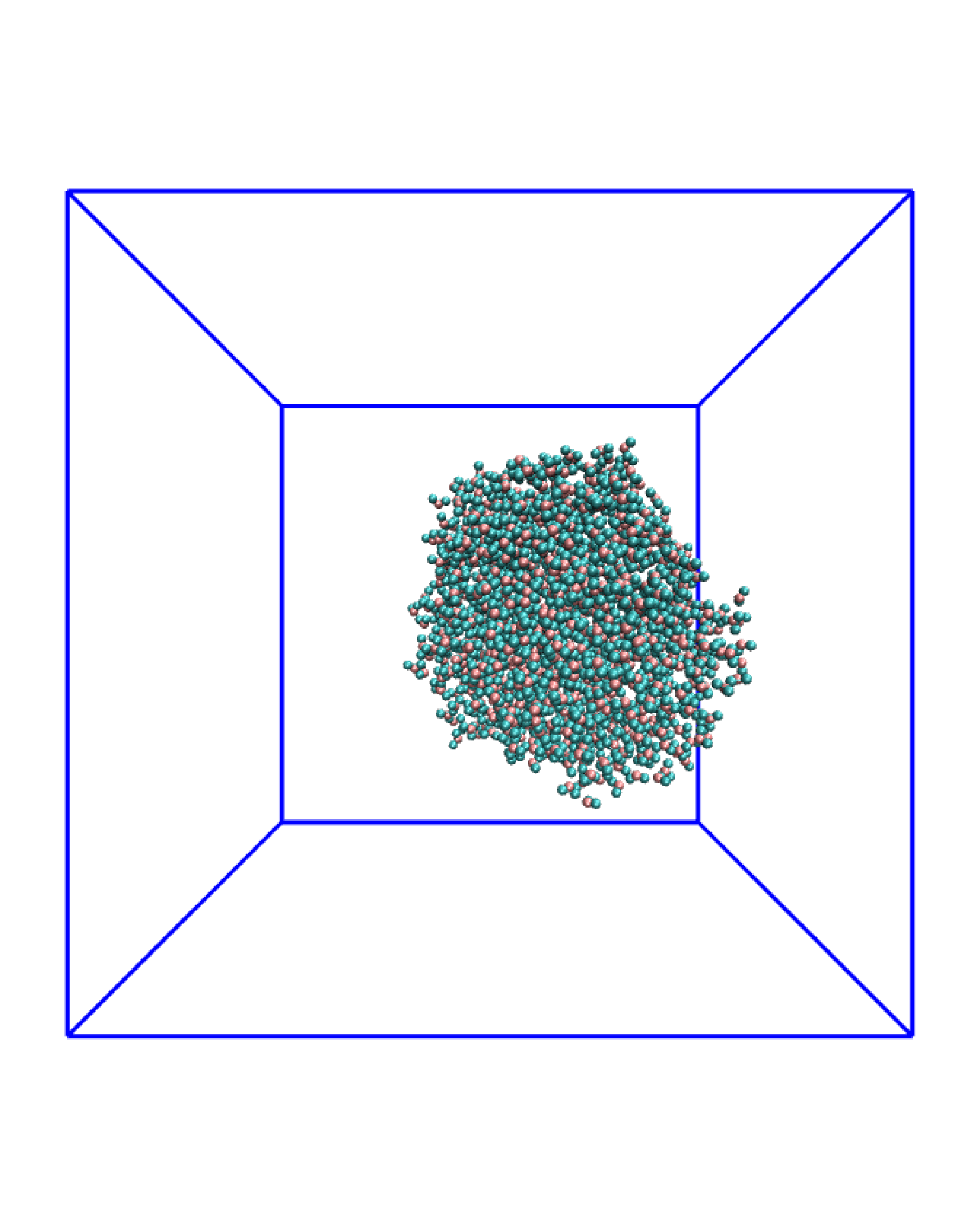}}\\
\caption{
(a) Initial artificial (low density) FCC crystal structure of H$_{2}$O system. (b) Disordered, liquid-like state of H$_{2}$O  system at $5\:ps$ simulation time. (c) The equilibrium state of a H$_{2}$O nanodroplet at $1\:ns$ simulation time.}

\label{fcc initial snapshot}
\end{figure}

We model the evolution of a molecular system of liquid water (H$_{2}$O), under constant volume, in a transient time regime, starting from a crystal  structure of low density towards the formation of a nanodroplet. The system consists of $1099$ H$_{2}$O molecules in a cubic box at temperature $T=300\:K$. The initial positions of atoms are at an FCC crystal structure; see Fig. \ref{fcc initial snapshot}(a). 
The FCC unit cell length is $10$ \AA. The box size is $70^{3}$ \AA$^{3}$, while the density of the system is $0.09577\: \frac{gr}{cm^{3}}$. We use the SPC/E (extended simple point) atomistic force field characterized by the parameters listed in table \ref{Force field parameters nvt} \cite{SPCE}. The SPC water model specifies a 3-site rigid water molecule with charges on each of the three atoms.

We perform AA MD simulations under the \textit{NVT} ensemble at temperature $T=300\:K\:$ using the \textit{Nose-Hoover} thermostat to keep the temperature constant. We should mention herein that the transition of the system out of equilibrium towards equilibrium depends on the choice of the statistical ensemble. Here, we perform simulations under the \textit{NVT} ensemble, since it is closer to a real experiment, where it is more convenient to control the temperature rather than the energy. Additionally, \textit{NVE} simulations are usually more sensitive with respect to numerical issues, than the \textit{NVT} ones. That is, the truncation errors during the integration process lead to fluctuation in the total energy, and more importantly to a final temperature that is not the desired one. The above-mentioned were confirmed by performing additional MD simulations in the \textit{NVE} ensemble, starting from the same initial conditions (data not shown here). We use the $\textit{velocity Verlet}$ time integrator with the time step $1\:fs$. The long-range Coulombic interactions were calculated using the PPPM solver. The SHAKE procedure was used to conserve intramolecular constraints \cite{29acd3d494044594aea0829ef236aad6,LAMMPS}.

 \begin{table}[htbp]
\centering
\resizebox{\textwidth}{!}{
\begin{tabular}{|c|c|c|c|c|c|c|c|c|c|}
\hline
   $\epsilon_{OO} \:[\frac{Kcal}{mol}] $& $\sigma_{OO}$ [\AA]  &$\sigma_{OH},\epsilon_{OH} $&$\sigma_{HH},\epsilon_{HH} $&bond length\: [\AA]&$k_{b}\:[\frac{kcal}{mol {\AA}^{2}}]$ &angle $\theta\: [deg] $&$k_{\theta}\: [\frac{kcal}{mol rad^{2}}]$&$q_{O}\:[e]$&$q_{H}\:[e]$\\
\hline
$0.1553$  & $3.16562$ &$ 0$&$0$&$1.0$&$412.00935$ & $109.47$ &$45.739$ &$-0.8476$& $+0.4238$\\
\hline
\end{tabular}}
\caption{Force field parameters of SPC/E H$_{2}$O system.}
\label{Force field parameters nvt}
\end{table}

The PSFM methods, the  radial pair distribution function (RDF) between H$_2$O molecules, and the self-diffusion coefficient (DC) of H$_2$O molecules were computed with an in-house simulation package using \textit{Matlab} \cite{MATLAB:2017}. The required MD simulations were performed using the $\textit{Lammps}$ package \cite{LAMMPS}.

We generate four atomistic data sets with characteristics listed in table \ref{data set table nvt}. The $1^{st}$ and the $2^{nd}$ data sets consist of paths corresponding to the system's transient regime. The initial configurations $\{q_{0,n}\}=\{q_{0}\}\in\mathbb{R}^{3N},n=1,..,n_{p}$ correspond to the same FCC crystal structure as discussed above, see Fig. \ref{fcc initial snapshot}(a), and  the initial velocities  were assigned according to the Boltzmann distribution corresponding to the requested temperature of $300\:K$. The $1^{st}$ 
data set has configurations corresponding to the time interval $[0,100]\: ps$. It is used to estimate the significant part of the transient regime empirically. The $2^{nd}$ data set corresponds to a shorter time interval $[0,10]\: ps$ and is the one that represents the significant part of the transient regime. The latter is used as a sample set to estimate the time-dependent CG interaction potential characterizing the proposed model \eqref{eq:Langevin_time_d} as well as to estimate the friction coefficient based on \eqref{eq: transient friction formula}. The $3^{rd}$ and $4^{th}$ data sets consist of paths corresponding to the equilibrated system. 
The $3^{rd}$ data set consists of independent and identically distributed configurations. It is used to estimate the equilibrium CG interaction potential and the corresponding RDF. In contrast, the $4^{th}$ one consists of correlated configurations and is used to estimate the equilibrium friction coefficient given by \eqref{eq: equilibrium friction formula}.

\begin{table}[htbp]
\centering
\resizebox{\textwidth}{!}{
\begin{tabular}{|c|c|c|c|c|c|c|}
\hline
Data sets& Molecules $M$& Paths $n_{p}$&Observation frequency [ps] &Observation time interval $\tr$ [ps]& Configurations per path $n_{t}$\\
\hline
$1^{st}$ data set   &$1099$& $20$ & $0.1$ & $[0,100]ps$ & $1000$ \\
\hline

$2^{nd}$ data set   &$1099$& $10$ & $0.1$ & $[0,10]ps$ & $100$ \\
\hline

$3^{rd}$ data set &$1099$& $10$ & $0.1$ & $[1000,1100]ps$ & $1000$ \\

\hline
$4^{rth}$ data set &$1099$& $1$ & $0.05$ & $[1000,1100]ps$ & $2000$ \\

\hline
\end{tabular}}
\caption{Details of the generated data sets from the all-atom MD simulations for the H$_{2}$O molecular system.}
\label{data set table nvt}
\end{table}

\subsection{Estimation of the short-time dynamics}\label{short-time dynamics subsection}
\subsubsection{Empirical estimation of the transient regime}
First, we estimate  the time interval corresponding to the significant part of the system's transition from a crystal to a liquid state by tracking the time-dependent CG pair potential and the atomistic RDFs.

 We solve the PSFM optimization problem given in Eq.~\eqref{PSFM instant times}, for several instantaneous times $t=0.5,2,5,10,30,50,100\:ps $, to obtain the pair interaction potential $\bar{u}(r;\phi^{\ast}_t)$. The input data set in the optimization problem is the $1^{st}$ data set in table \ref{data set table nvt}. For the parametrization of the CG pair potential with respect to the distance $r$, we use $N_{d}=48$ linear splines.
 
 In Fig. \ref{evolution potential FM, nvt water}, we present the time evolution of the effective CG potential, obtained by solving the above-mentioned optimization problem. In the same graph, we  plot the equilibrium CG pair  interaction potential (approximating the PMF), $\bar{u}(r;\phi^{\ast})$, calculated by solving the PSFM method using the $3^{rd}$ data set in table \ref{data set table nvt}. 
We observe that as time progresses the minimum in the interaction potential changes with time while its location remains constant. 
This result depicts that as time increases 
the system gradually reaches equilibrium and the derived CG pair potential approaches the PMF.
The strong attraction between molecules (large values of the minimum in the potential) at short times indicates that the molecules start from an initial condition where all interactions (and forces) are attractive; it also exhibits the tendency of the H$_2$O molecules to attract each other in order to form the equilibrium (nanodroplet) state under the specific thermodynamic conditions. 

 \begin{figure}[htbp]
 \centering
\includegraphics[width=0.8\textwidth]{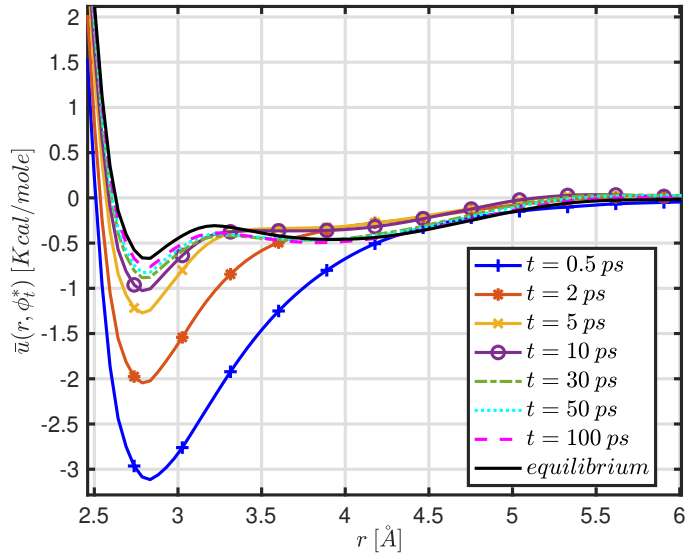}
 \caption{Time evolution of the effective CG pair potential representing the system free energy. As time increases, the effective CG potential reaches the equilibrium potential.}
    \label{evolution potential FM, nvt water}
\end{figure}

In addition to the time evolution of the CG pair potential, we estimate the time evolution of the atomistic RDF analyzed at the level of the center of mass of the CG particles, see Fig. \ref{evolution g(r), nvt water}. 
We observe that (a) at $t=0\:ps$ the local density equals the average density. This result confirms that the molecules are initially in a sparse condition, occupying the entire space of the simulation box. It is clear also that as time increases the transition towards liquid takes place and the RDF becomes greater at the distance of $[2.5,12.5]\:\AA$. For large instant times, i.e., at $t>30\ ps$ and at short pair distance $r$ the RDF increases, while at a long distance, the RDF converges to zero. This is because the strong cohesive forces between molecules result in their coming closer, creating a nanodroplet, see Fig. \ref{fcc initial snapshot} (c) and Fig. \ref{evolution g(r), nvt water} (b).

\begin{figure}[htbp]
\centering
\subcaptionbox{}[.49\textwidth]{\includegraphics[width=0.49\textwidth]{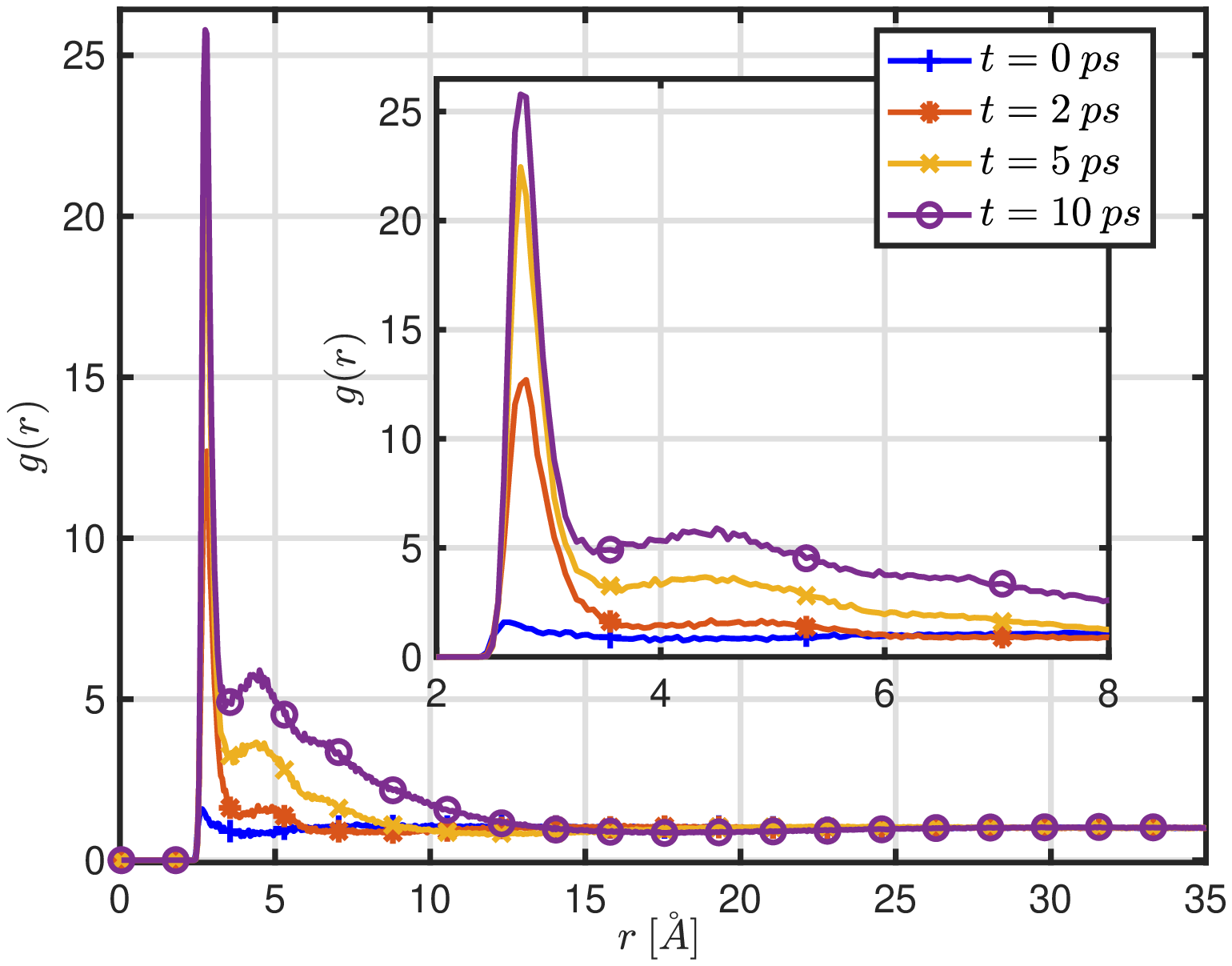}}
\subcaptionbox{}[.49\textwidth]{\includegraphics[width=0.49\textwidth]{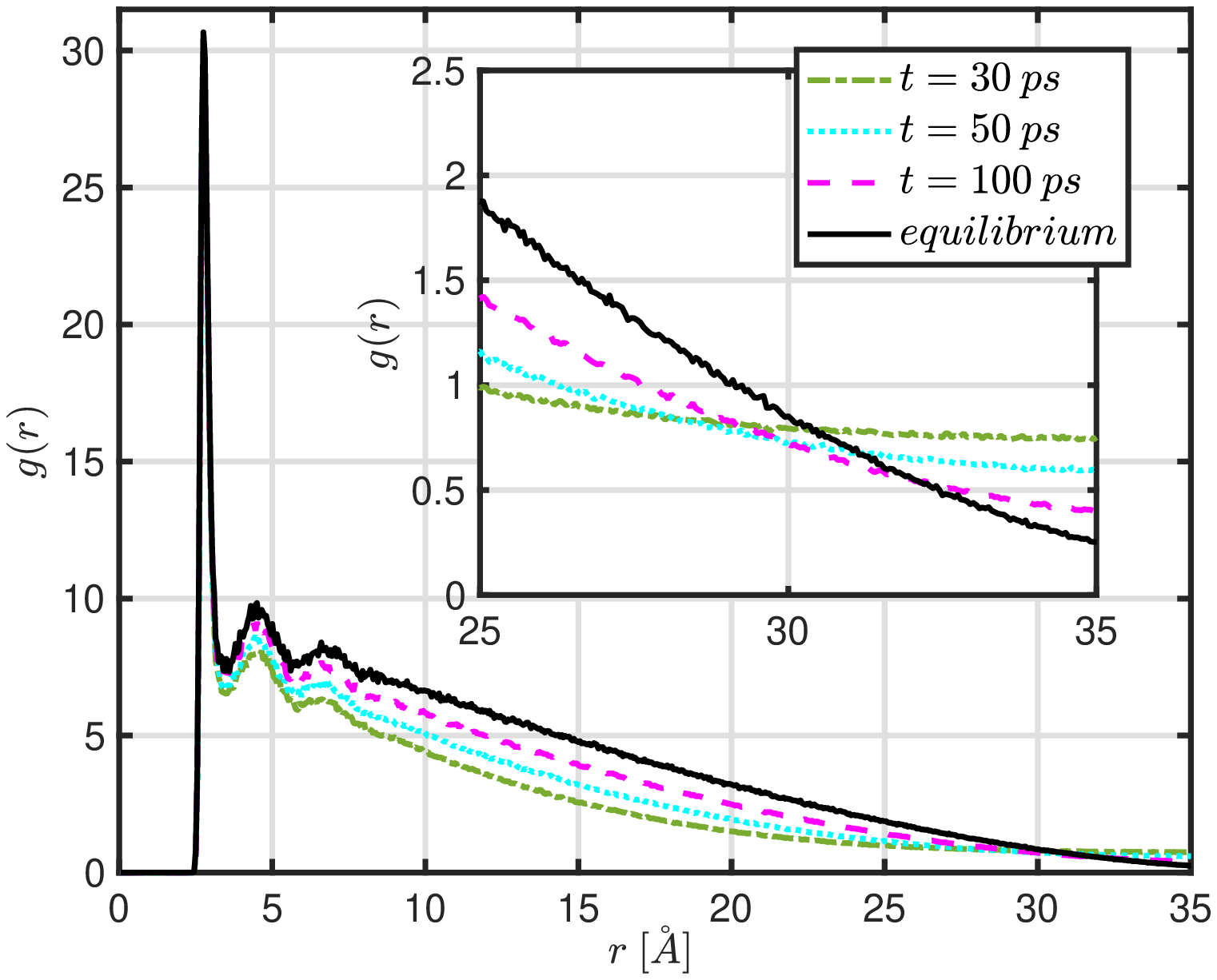}}
\caption{Time evolution of the atomistic RDF $g(r)$ at instantaneous times. (a) Early times (b)  Longer times. In the inset is the enlarged graph. As time increases, the molecules come closer creating a nanodroplet and the $g(r)$ converges to the equilibrium one.
}\label{evolution g(r), nvt water}
\end{figure}
 
Based on these results, see Fig. \ref{evolution potential FM, nvt water} and \ref{evolution g(r), nvt water}, we assume that the significant part of the transient regime of the H$_{2}$O system,  $\tr$, corresponds to the time interval of $[0,10]\:ps$. Hereafter, we produce the $2_{nd}$ data set in table \ref{data set table nvt}, which is utilized to estimate both the time-dependent CG potential and the friction coefficient.

\subsubsection{Estimation of the time-dependent coarse-grained interaction potential}

We estimate the time-dependent CG interaction potential, $\bar{u}(r,t;\theta^{\ast}_{\tr})$, by solving the PSFM problem, given in \eqref{splines exp. on t and r, opt.problem}, using four paths of the $2^{nd}$ data set. For the parametrization of the time-dependent CG potential with respect to time $t$ and distance $r$, we used linear splines with $N_{b}=24$ and $N_{d}=48$, respectively. Thus, the total number of parameters in the optimization problem was $24\times 48=1152.$ Figure \ref{evolution potential PSFM, nvt water} depicts the time evolution of the effective time-dependent CG pair potential. As time increases  the derived CG pair potential approaches the PMF, similarly to the findings described in Fig. \ref{evolution potential FM, nvt water}.

 \begin{figure}[htbp]
 \centering
\includegraphics[width=0.8\textwidth]{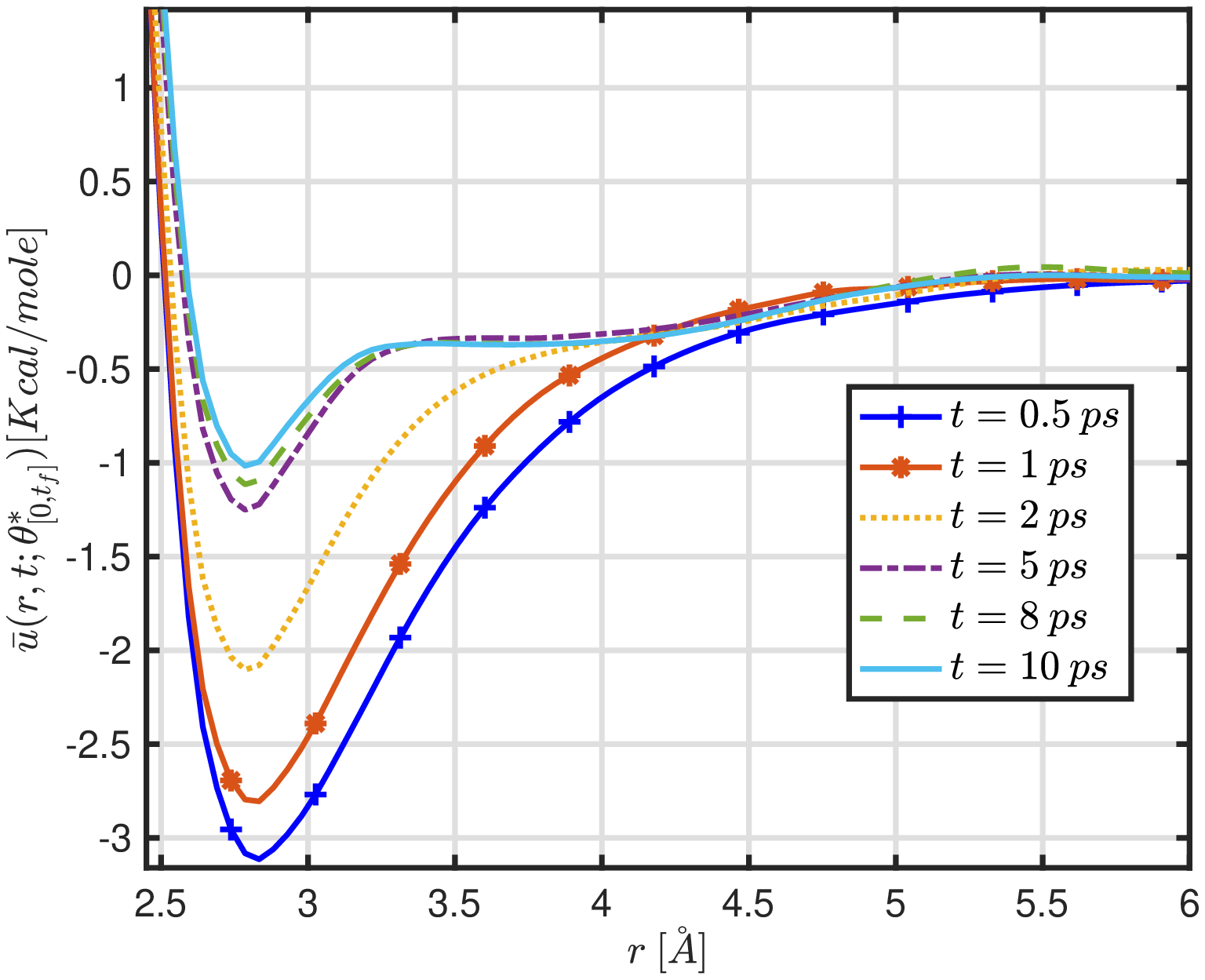}
  \caption{Time evolution of the effective time-dependent CG pair potential derived through the PSFM method given in \eqref{splines exp. on t and r, opt.problem}. }
    \label{evolution potential PSFM, nvt water}
\end{figure}

\begin{figure}[htbp]
 \centering
\includegraphics[width=0.8\textwidth]{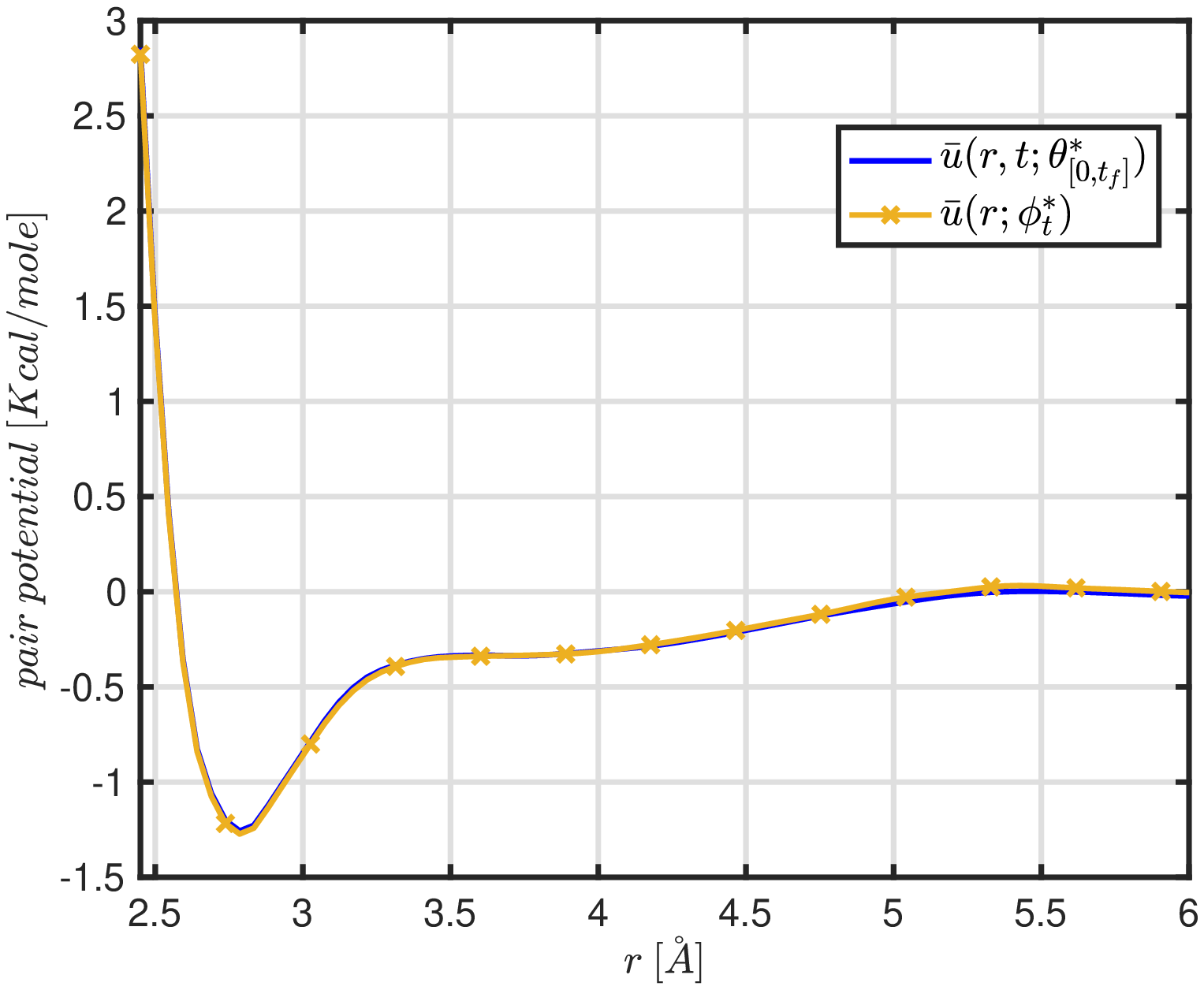}
 \caption{Comparison at $5\: ps$ of the effective pair potentials as derived by solving the PSFM method for fixed instantaneous time and parametrizing with respect to the pair distance $r$ and by solving the  PSFM problem in which both the time $t$ and the pair distance $r$ are variables. We get similar results for other instantaneous times as well.} \label{compare eff potential at 5 ps nvt water}
\end{figure}

 Fig. \ref{compare eff potential at 5 ps nvt water} indicates that if one solves the PSFM method for fixed instantaneous time, parameterizing the CG pair interaction potential with respect to the pair distance $r$, and the PSFM method, in which the CG potential is parameterized with respect to the pair distance and time, the effective CG pair potentials are equivalent. 
The difference between the two optimization methods lies in the fact that for solving the PSFM optimization problem for a fixed instantaneous time, e.g., at $5\:ps$, we use as a sample the configurations corresponding to the specific instantaneous time extracted from twenty paths, while for the case we represent the CG potential as a time-dependent quantity we use configurations from four paths at $[0,10]\:ps$. Despite the high dimensionality of the latter, it is more advantageous as it provides the time evolution of the effective pair potential only by solving the optimization problem once. In contrast, the first optimization problem requires the solution of the minimization problem multiple times for the derivation of the pair potential's time evolution.
We verified this equivalence by comparing the CG potential at several instantaneous times of the transient regime. 

\subsubsection{Estimation of the friction coefficient}
Next, we estimate the friction coefficient in two ways using two different data sets. In the first case, we estimate the friction coefficient empirically based on the formula \eqref{eq: transient friction formula}
using the $2^{nd}$ data set. The value of the friction coefficient we derived is
$\ft=0.0059\frac{gr}{mole}\frac{1}{fs}.$ In the second case, based on \eqref{eq: equilibrium friction formula} to calculate the time correlation functions for approximating the friction kernel, we use the $4^{rth}$ data set that corresponds to equilibrium. The value of the friction coefficient we derived is $ \fe=0.4484\frac{gr}{mole}\frac{1}{fs}$, which is similar to the diffusion coefficient value of the bulk water at $\mathrm{P}=1\:atm$. That is the interactions between the atoms in the nanodroplet resemble their interactions in the bulk water.
As expected, we observe a deviation in the value of the estimated friction coefficient because the molecular system under study is non-homogeneous. The friction that arises from collisions between atoms is
less intense when they are initially in a sparse condition due to the FCC structure, while the friction is more vigorous when they come close, approaching the behavior of a realistic water nanodroplet.

\subsection{Model comparison and validation of the proposed CG model}
To strengthen the validity of the proposed CG model we compare the CG dynamics via  various models, given in table \ref{cg models}. 
The dynamics are  the  Langevin dynamics \eqref{eq:Langevin_time_d} and the models differ due to (a) the pair potential: the potential of mean force (\textit{PMF}) or the time-dependent (\textit{TD}) one, and (b) the friction coefficient: estimated  from the equilibrium (\textit{fe}),  the transient regime (\textit{ft}), or without friction (\textit{0}), see table \ref{cg models}.

For the \textit{TD-fe} model,  we perform CG MD simulations using as interaction potential the time-dependent CG potential, $\bar{u}(r,t;\theta^{\ast}_{[0,t_{f}]})$, and the  $\fe$  friction coefficient using equilibrium AA information.
Similarly, for the \textit{PMF-fe} model, we have the equilibrium CG pair potential $\bar{u}(r;\phi^{\ast})$ and the  $\fe$ friction coefficient.  
For the \textit{PMF-ft} and \textit{TD-ft} models  the friction coefficient $\ft$ is used, and the CG pairs potentials $\bar{u}(r;\phi^{\ast})$ and $\bar{u}(r,t;\theta^{\ast}_{[0,t_{f}]})$, respectively.
In the last two CG models, \textit{TD-0} and \textit{PMF-0}, the time-dependent CG interaction potential and the equilibrium CG potential are used, respectively, without friction in the dynamics. The MD simulations for these CG dynamics models were performed under the same conditions of the atomistic ones, i.e., under the \textit{NVT} ensemble at $300\:K$. Details about these CG MD runs are listed in table \ref{CG RUNS data set table nvt}.

We compare and validate the models with respect to two observables: the RDF which provides structural information, and the DC, which is related to the system kinetics.
To validate each model we compare it to the AA corresponding observables, while we define a better model to be the one that is "closer" to the AA corresponding observable.

\begin{table}[htbp]
\centering
\begin{tabular}{|c|c|c|}
\hline
CG Model &
CG interaction potential& Friction\\
\hline
\textit{TD-fe}  &
time-dependent potential, $\bar{u}(r,t;\theta^{\ast}_{[0,t_{f}]})$&$\fe$\\
\hline
\textit{PMF-fe} &
PMF, $\bar{u}(r;\phi^{\ast})$&$\fe$  \\

\hline
\textit{PMF-ft}  &
PMF, $\bar{u}(r;\phi^{\ast})$&$\ft$\\
\hline
\textit{TD-ft } &
time-dependent potential, $\bar{u}(r,t;\theta^{\ast}_{[0,t_{f}]})$&$\ft$\\
\hline
\textit{TD-0} &
 time-dependent potential, $\bar{u}(r,t;\theta^{\ast}_{[0,t_{f}]})$&-\\
\hline
\textit{PMF-0}  &
PMF, $\bar{u}(r;\phi^{\ast})$&-\\
\hline
\end{tabular}
\caption{CG dynamics models proposed and studied for the transient dynamics of H$_{2}$O molecular system.}
\label{cg models}
\end{table}

\begin{table}[htpb]
\centering
\resizebox{\textwidth}{!}{
\begin{tabular}{|c|c|c|c|c|c|c|}
\hline
Molecules $M$& Paths $n_{p}$&Observation frequency [ps] &Observation  time interval, $\tr$ [ps]& Configurations per path $n_{t}$\\
\hline
$1099$& $10$ & $0.1$ & $[0,10]$ & $100$ \\

\hline
\end{tabular}}
\caption{Details of the generated data sets of the CG MD simulations for the H$_{2}$O molecular system.}
\label{CG RUNS data set table nvt}
\end{table}





\begin{figure}[htbp]
\centering
\subcaptionbox{}[.49\textwidth]{\includegraphics[width=0.49\textwidth]{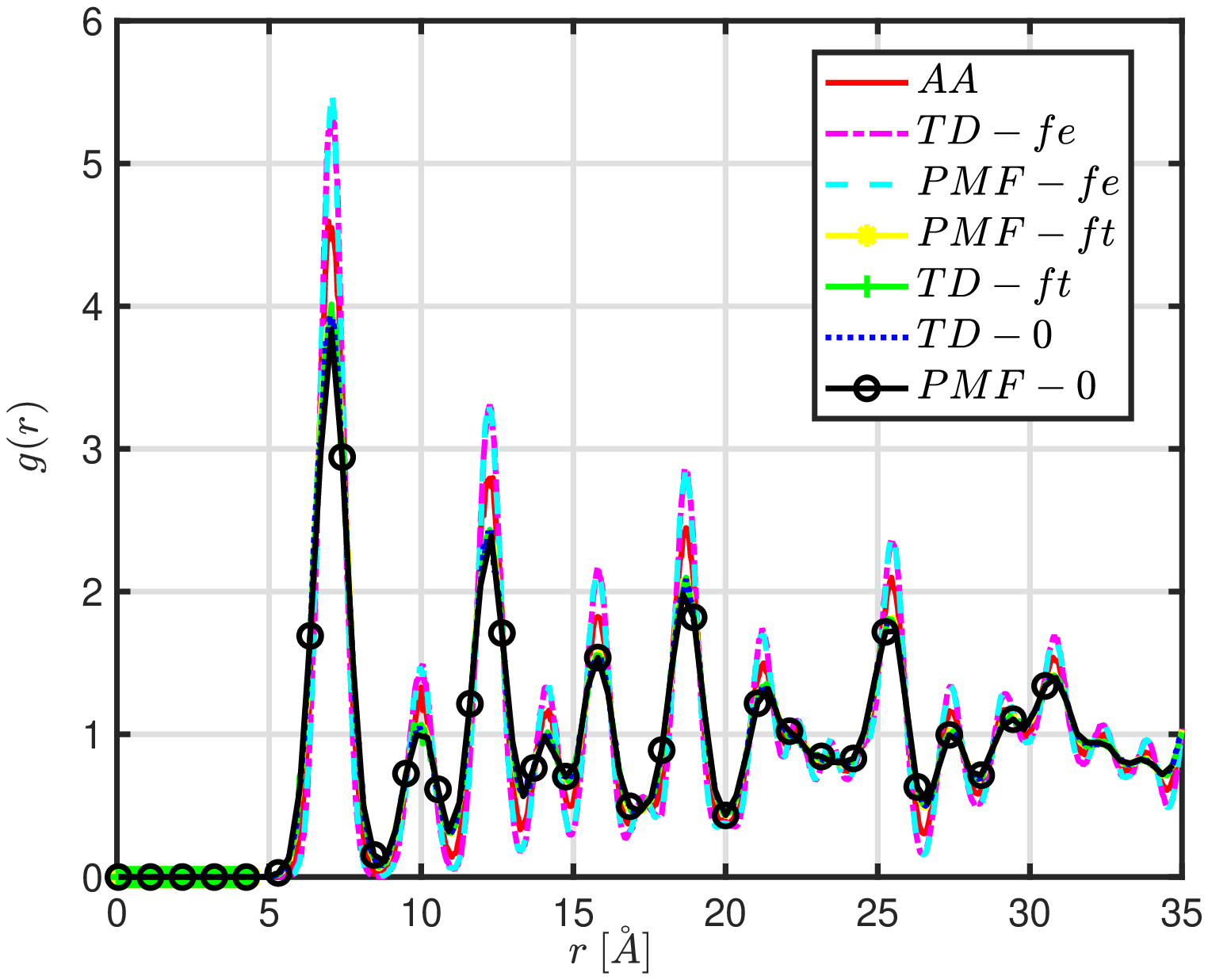}}
\subcaptionbox{}[.49\textwidth]{\includegraphics[width=0.49\textwidth]{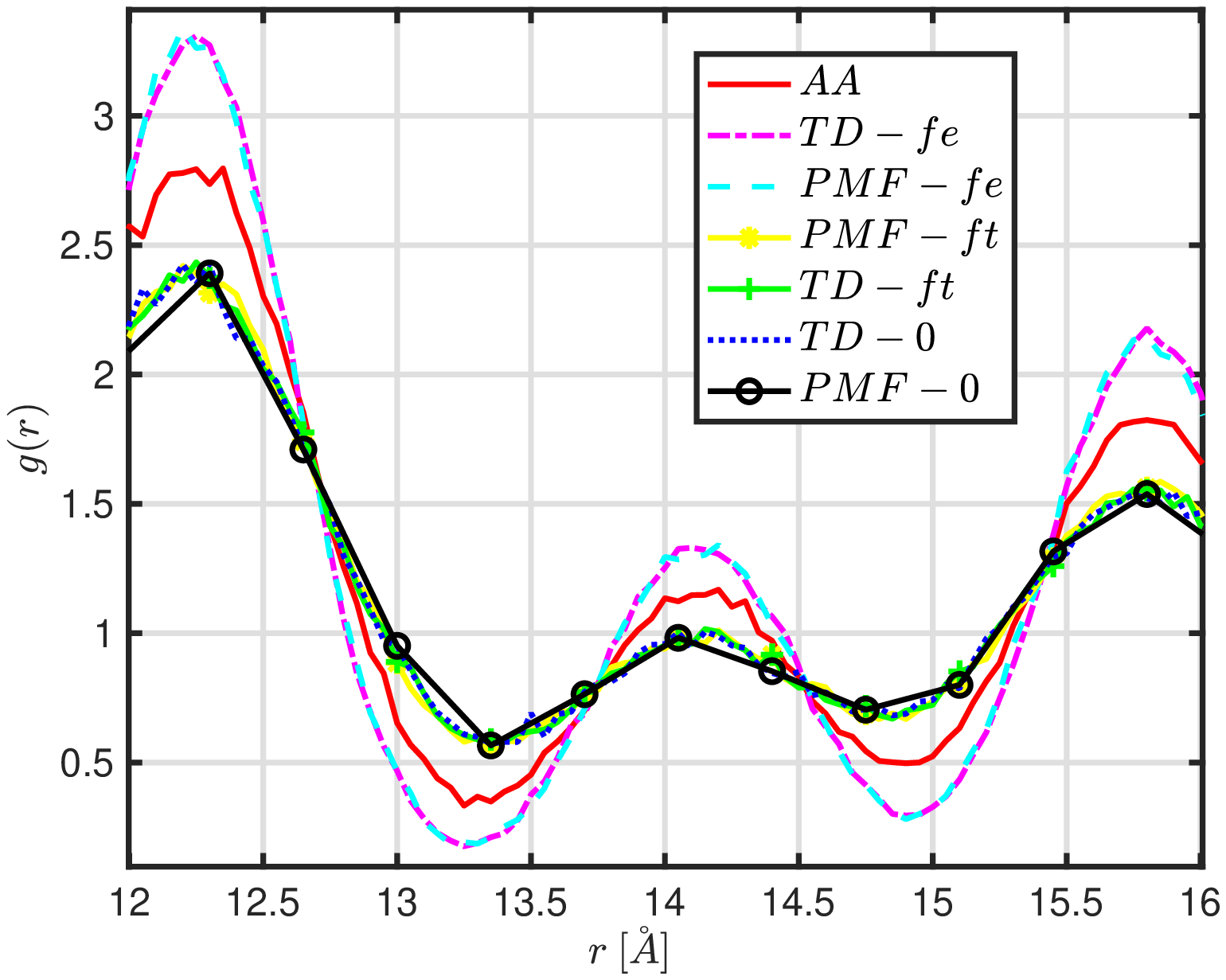}}
\caption{(a) RDFs of the water system at $0.1\:ps$, obtained from the AA and the CG simulations. (b) Enlarged graph. All CG models describe correctly the position of the peaks but not the corresponding amplitude.
}\label{compare gor at 0.1 allatom and cg models}
\end{figure}


\begin{figure}[htbp]
\centering
\subcaptionbox{}[.49\textwidth]{\includegraphics[width=0.49\textwidth]{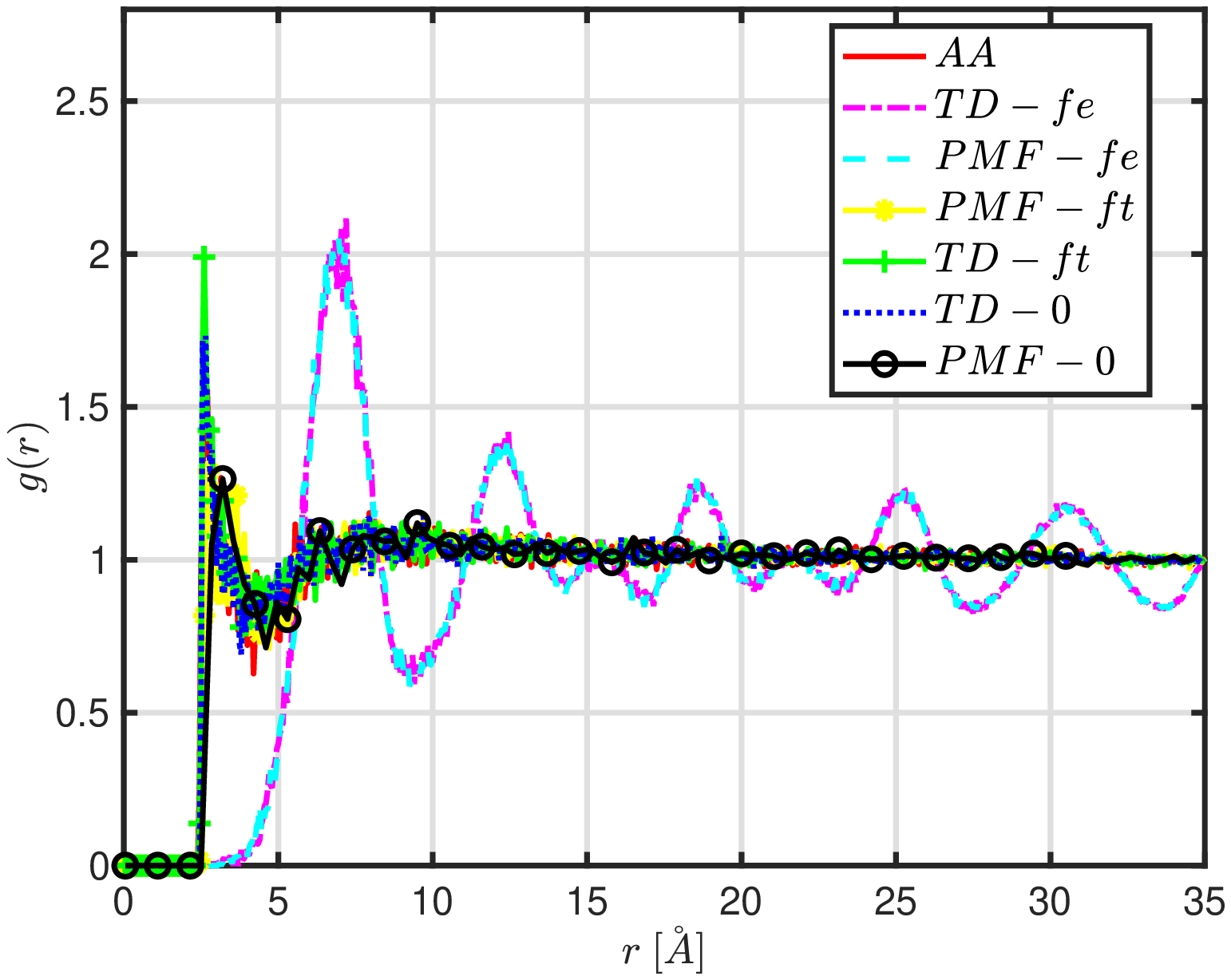}}
\subcaptionbox{}[.49\textwidth]{\includegraphics[width=0.49\textwidth]{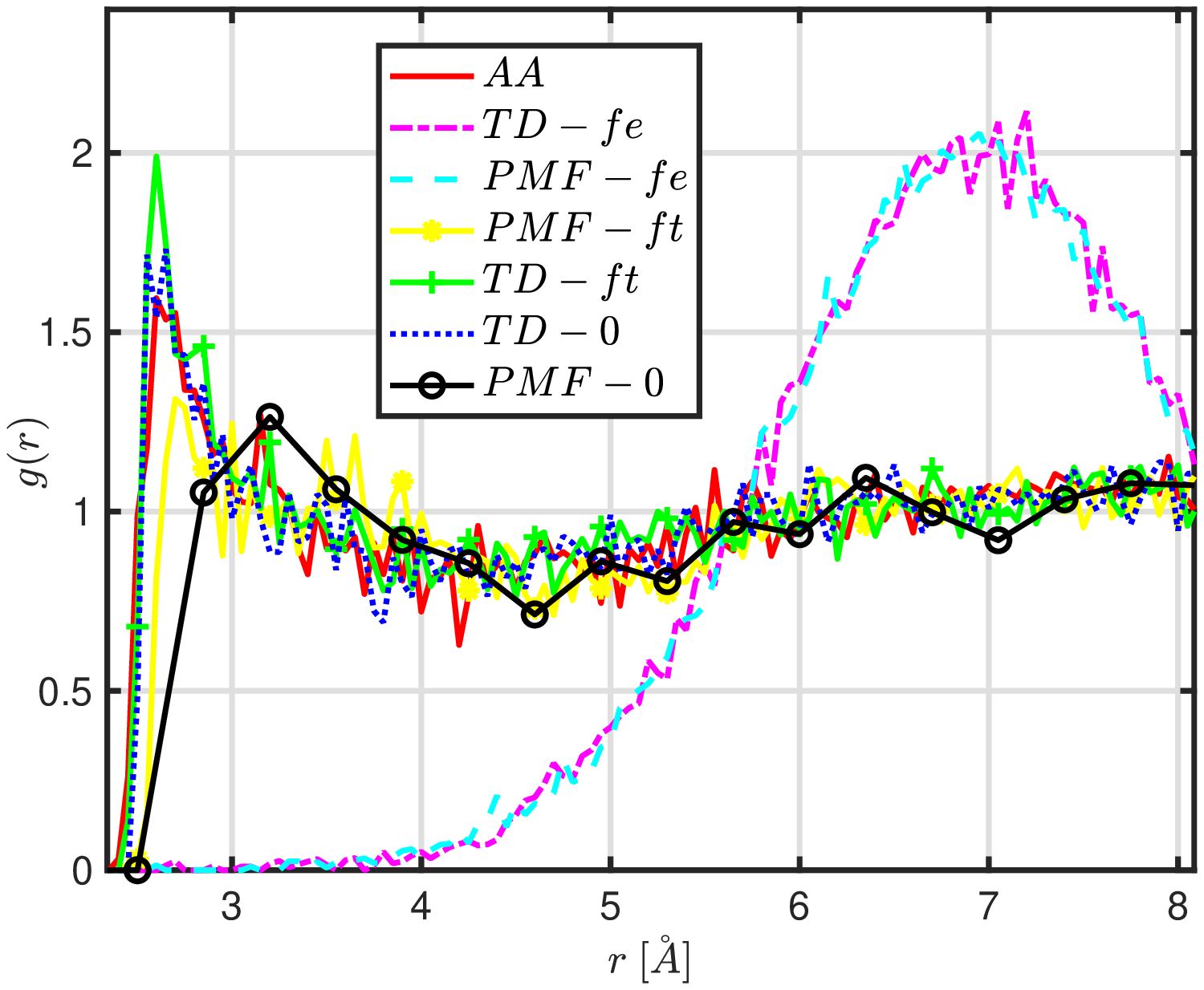}}
\caption{(a) RDFs of the water system at $0.5\:ps$, obtained from the AA and the CG simulations. (b) Enlarged graph. All models, except \textit{TD-fe} and \textit{PMF-fe}, adequately represent the system's transition towards the liquid state with slight deviations from the AA.
}\label{compare gor at 0.5 allatom and cg models}
\end{figure}

 \begin{figure}[htbp]
\centering
\subcaptionbox{}[.49\textwidth]{\includegraphics[width=0.49\textwidth]{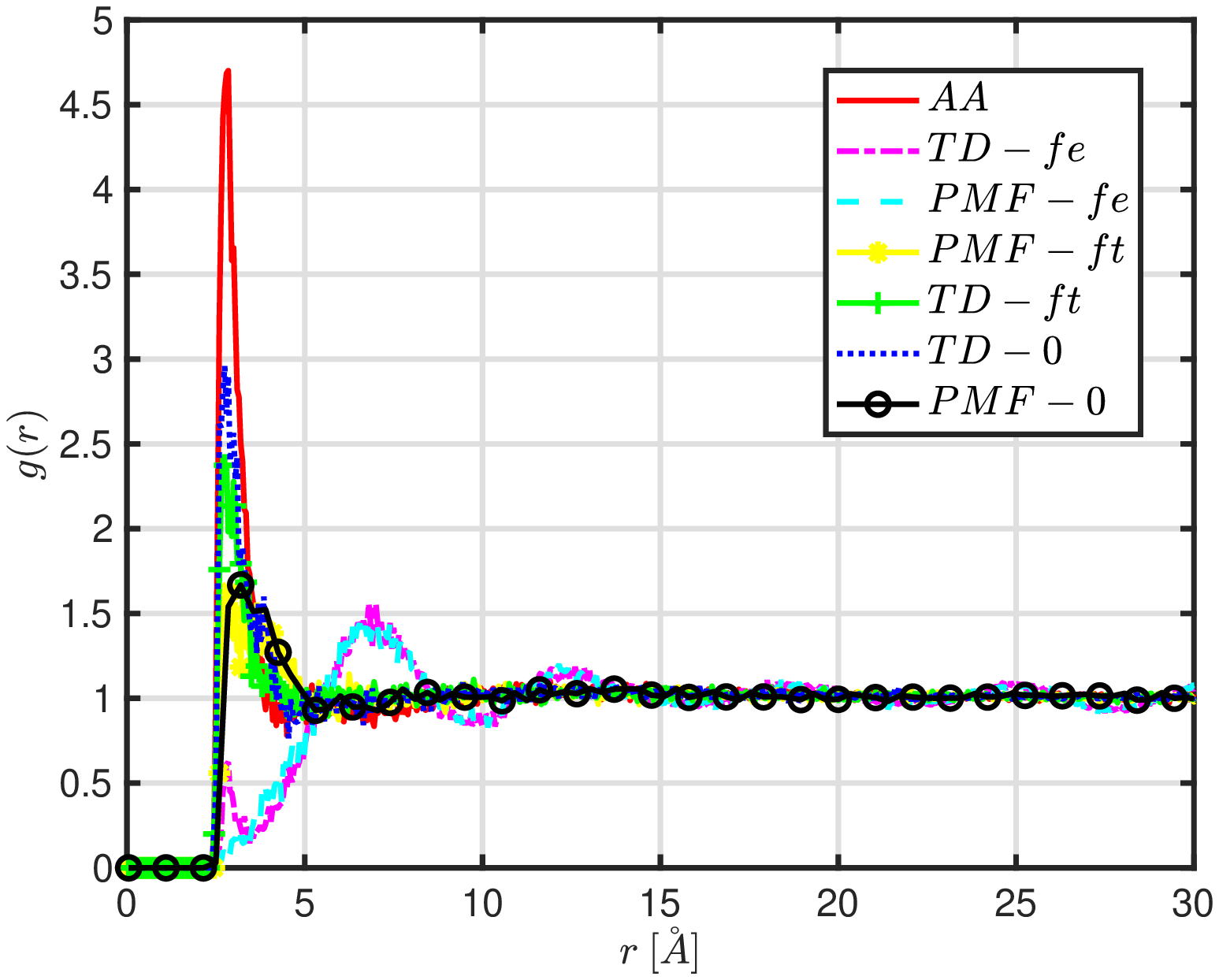}}
\subcaptionbox{}[.49\textwidth]{\includegraphics[width=0.49\textwidth]{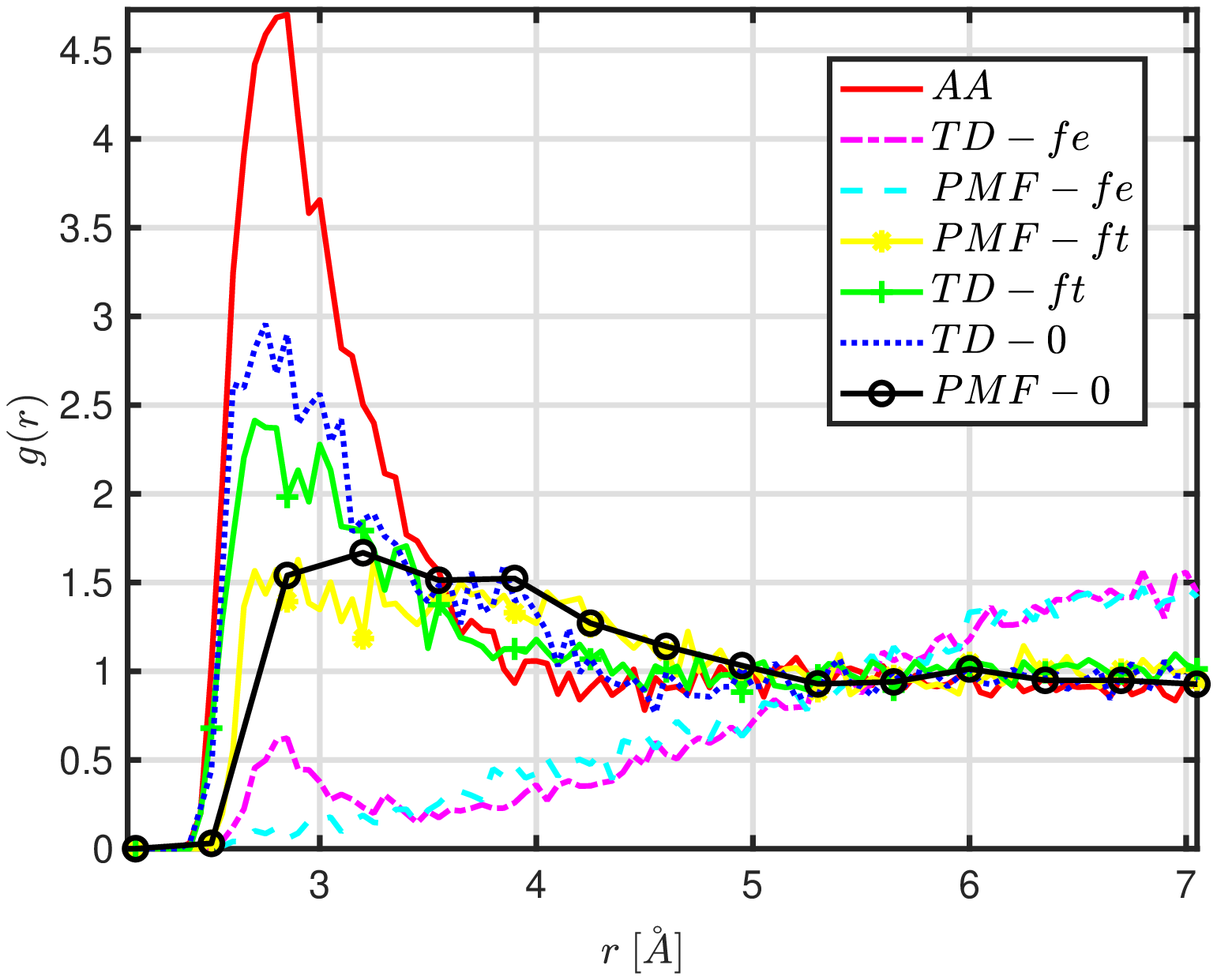}}
\caption{(a) RDFs of the water system at $1\:ps$, obtained from the AA and the CG simulations. (b) Enlarged graph. CG models \textit{TD-ft} and \textit{TD-0} are more representative approximations of the AA RDF relative to the other CG models as the corresponding peaks of the RDFs are closer to the AA one, and they still exhibit good agreement with this for large distances, above $5$\AA.
}\label{compare gor at 1 allatom and cg models}
\end{figure}
 
 \begin{figure}[htbp]
\centering
\subcaptionbox{}[.49\textwidth]{\includegraphics[width=0.49\textwidth]{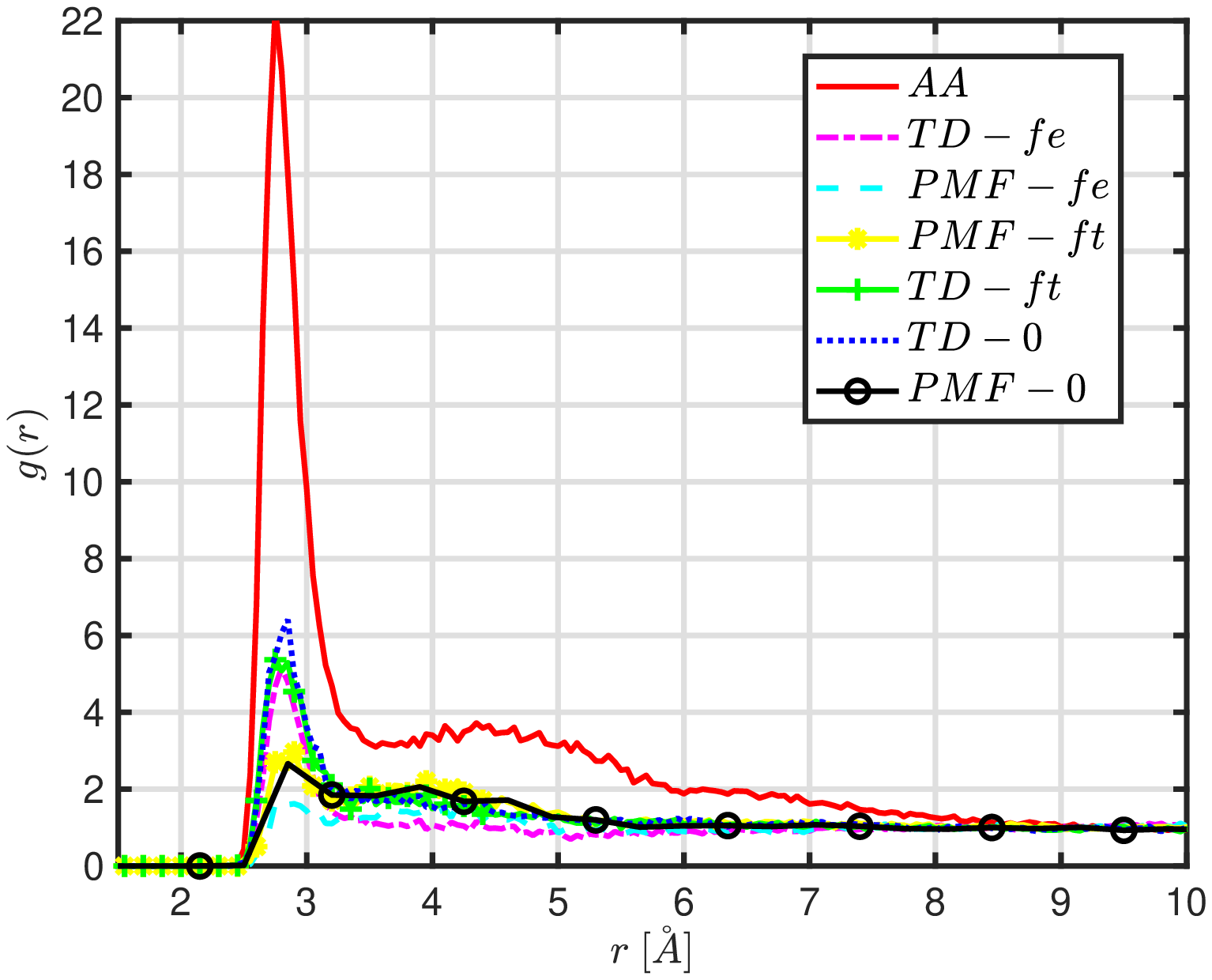}}
\subcaptionbox{}[.49\textwidth]{\includegraphics[width=0.49\textwidth]{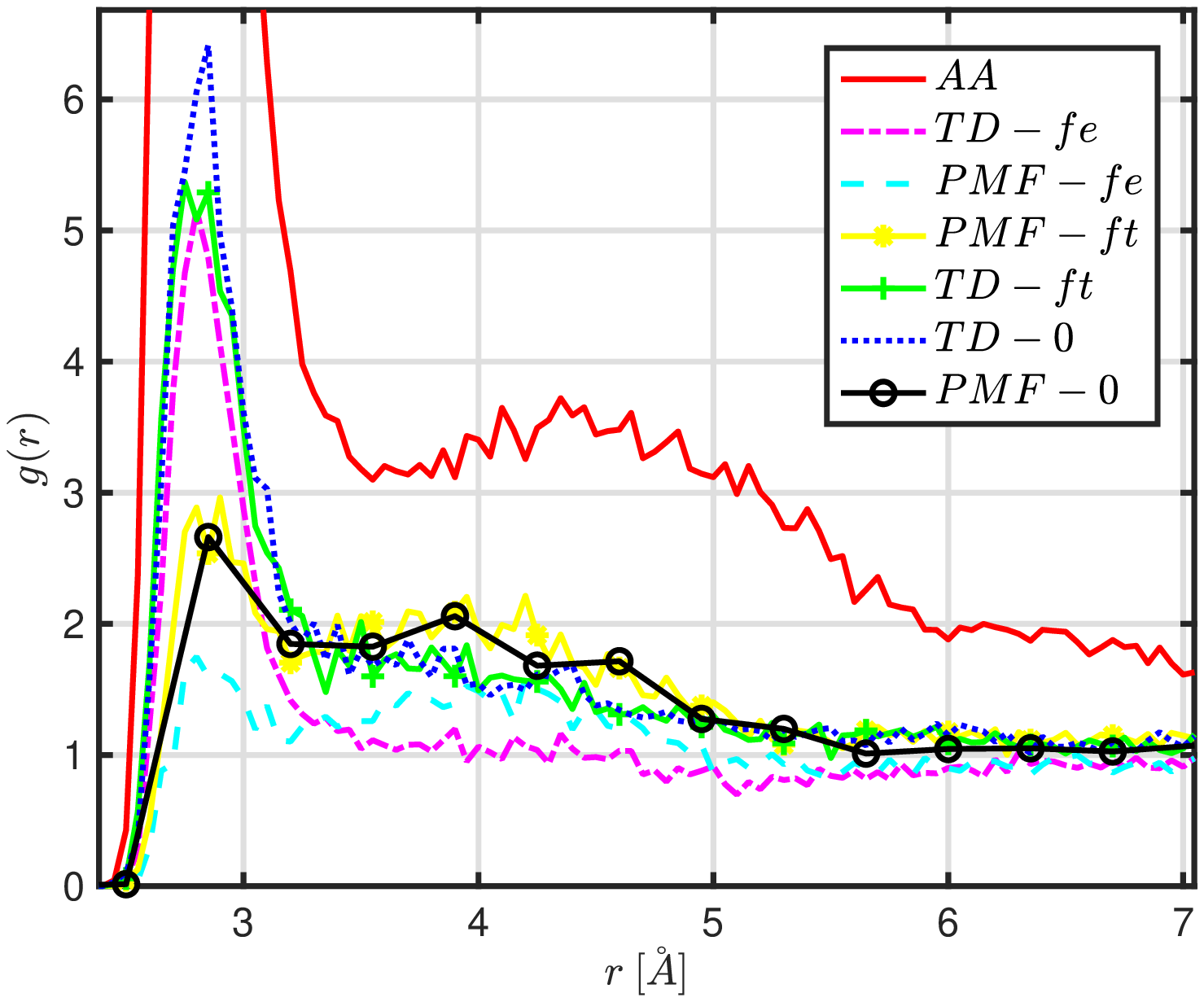}}
\caption{(a) RDFs of the water system at $5\:ps$, obtained from the AA and the CG simulations. (b) Enlarged graph. There are apparent deviations between the RDFs of the CG models and the AA one in all distances, especially at the position of the first peak. However, the \textit{TD-ft} and \textit{TD-0} models, which produce quite similar RDFs, are the most effective. There is also a significant improvement of the \textit{TD-fe} model's RDF representation.
}\label{compare gor at 5 allatom and cg models}
\end{figure}
 \begin{figure}[htbp]
\centering
\subcaptionbox{}[.49\textwidth]{\includegraphics[width=0.49\textwidth]{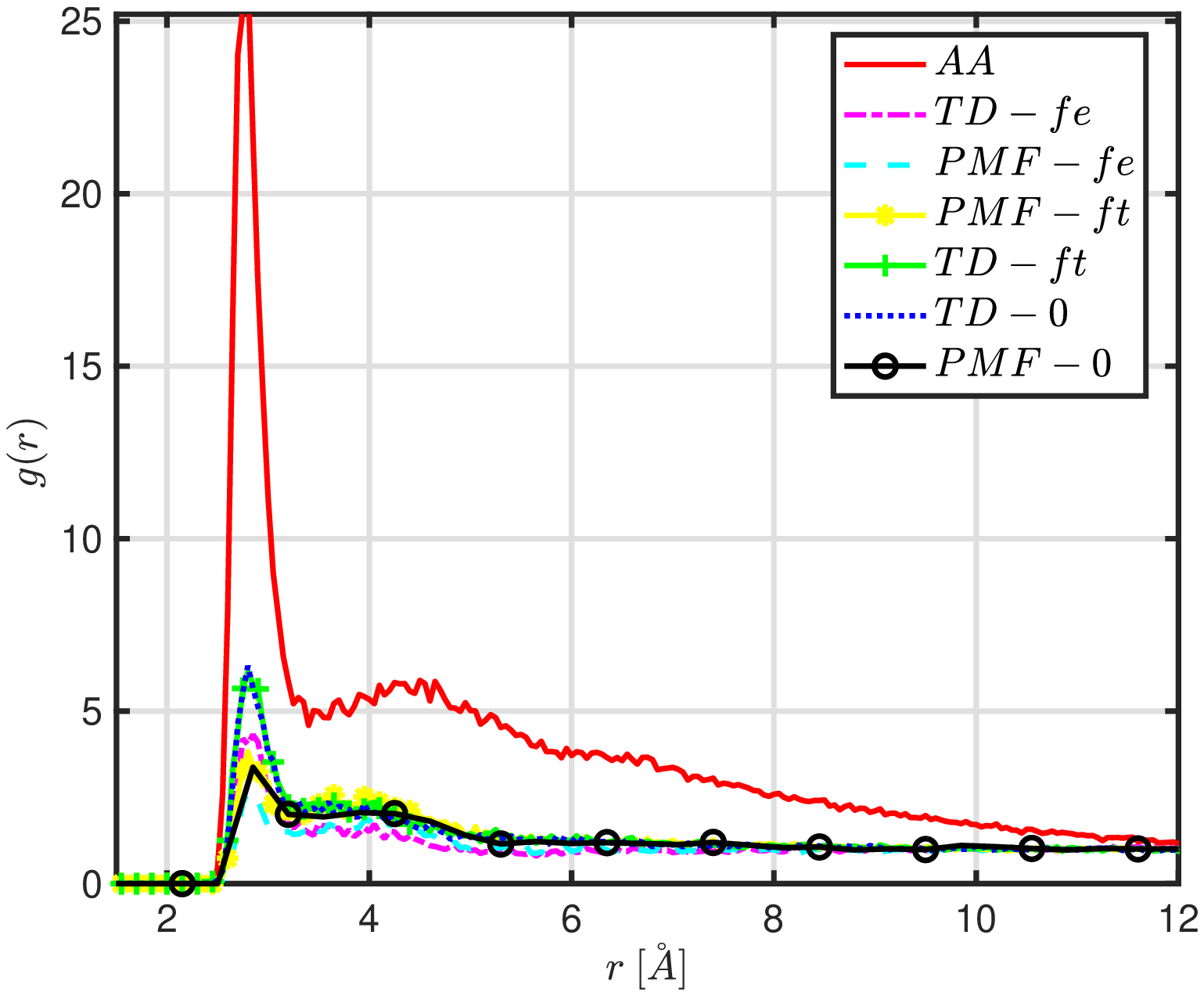}}
\subcaptionbox{}[.49\textwidth]{\includegraphics[width=0.49\textwidth]{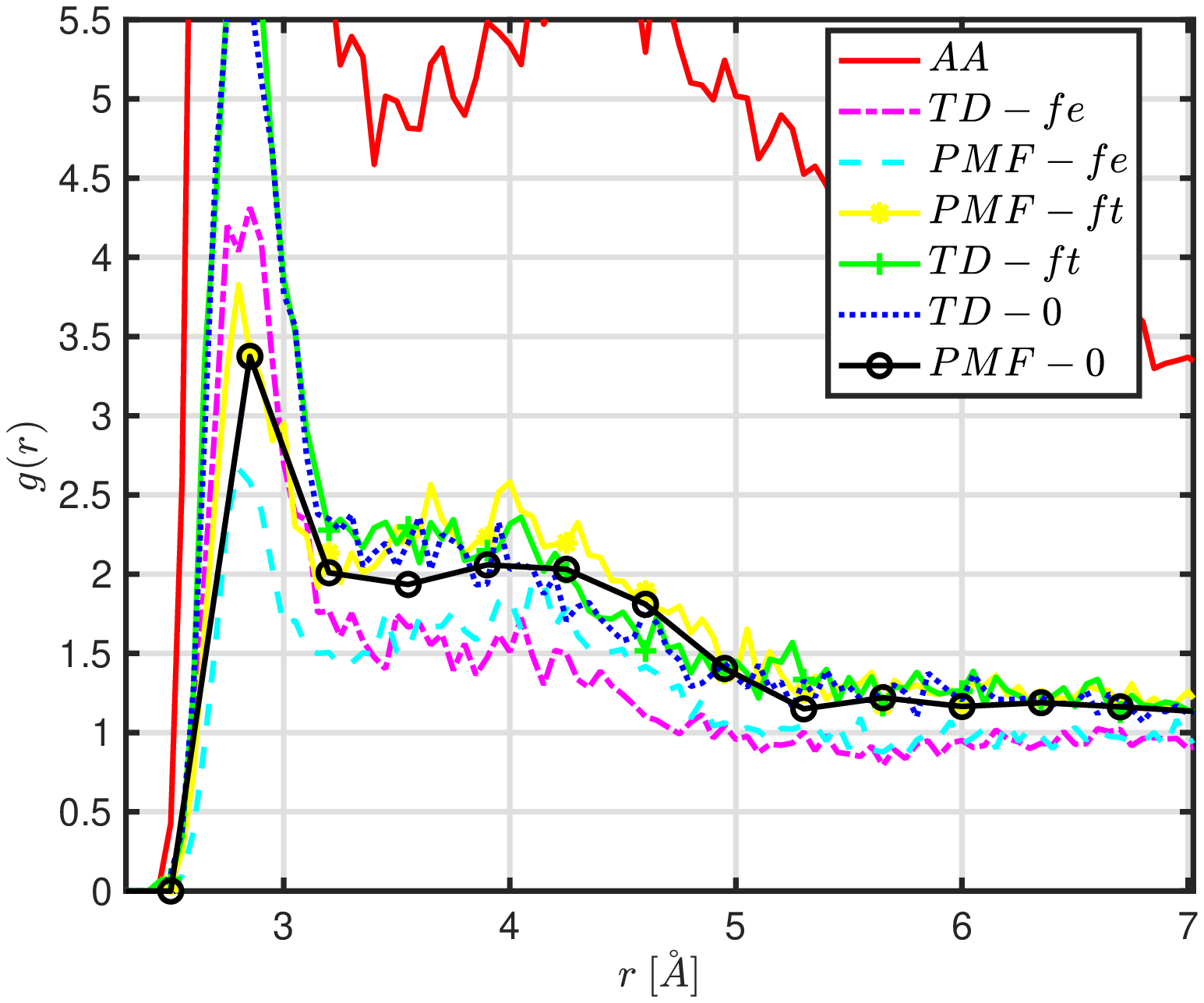}}
\caption{(a) RDFs of the water system at $10\:ps$, obtained from the AA and the CG simulations. (b) Enlarged graph. \textit{TD-ft} and \textit{TD-0}, whose RDFs are almost identical, are slightly more representative compared to the other CG models. However, the deviation of \textit{TD-ft} model's RDF from the AA RDF is larger for all distances compared to the corresponding one at shorter instant times. }\label{compare gor at 10 allatom and cg models}
\end{figure}

Figures \ref{compare gor at 0.1 allatom and cg models}-\ref{compare gor at 10 allatom and cg models} display the RDFs curves obtained from all the CG models reported in table~\ref{cg models} and their comparison with the corresponding AA one, at several \textit{instant times} of the transient regime. The AA RDFs are analyzed at the level of the center of mass of the CG particles. For the calculations of the RDFs at every instant time, we averaged over ten paths.
At very early instantaneous times of the transient regime, e.g., at $0.1,\:0.5\:ps$, the  \textit{PMF-0}   model seems to approximate well the AA RDF.  As time increases, the  \textit{PMF-0} model's effectiveness is weakened; this is not surprising since, due to the cohesive forces, water molecules, as time passes, approach each other to form a nanodroplet and the friction forces in the molecular system becomes more intense. 
Next, we observe that the models with the $\fe$ friction from equilibrium,  \textit{TD-fe} and \textit{PMF-fe}, generally result in an ineffective estimation of the system's short-time transient dynamics. Specifically, at $0.1\:ps$, both CG models produce RDFs with the same peak positions as the target AA RDF  but differ in amplitude size. At $0.5\:ps$, the shapes of the RDFs are entirely different, i.e., the curves obtained from the \textit{TD-fe} and \textit{PMF-fe} models present more peaks than the AA one. Nevertheless, the \textit{TD-fe} model is more effective than the \textit{PMF-fe} since, at $1\:ps, 5\:ps$ and $10\:ps$, it produces RDFs with the exact location of the first peak and size of amplitude closer to the target one. 
This is a clear indication that the widely used Langevin model given in \eqref{eq:Langevin equation} fails to represent the system's transient dynamics effectively.
Moreover, by comparing the  \textit{PMF-fe} and \textit{PMF-ft} models it is clear that if one uses the equilibrium interaction potential, estimating the friction coefficient using the AA data at the transient regime improves the dynamics. Therefore, such "transient-regime" friction is necessary if one is interested in CG transient dynamics.
We also observe that the \textit{TD-ft} model is consistently closer to the AA at all instant times than the \textit{PMF-ft} model. This result indicates the need for the \textit{time-dependent} pair interaction potential in representing the transient dynamics.
The comparison of the \textit{TD-ft} and \textit{TD-0} with the \textit{TD-fe} model demonstrates that the time-dependent CG potential in conjunction with the friction coefficient estimated using transient short-time AA data is the most efficient model among all the Langevin-type studied models for the representation of the short-time CG dynamics.
Furthermore, we observe that the RDF corresponding to the \textit{TD-0} model is very close to the one of \textit{TD-ft} model. In fact, at $0.5\:ps$, $1\:ps$ and $5\:ps$, the \textit{TD-0} model is slightly more effective. This observation   suggests that a better estimation of the friction coefficient used in the dynamics is necessary. 

At this point, it is worth mentioning that the effectiveness of the suggested model in representing the short-time dynamics of the system consisting of water molecules is not significantly affected if one performs simulations using the DPD algorithm, as we demonstrated by performing additional DPD simulations in a few specific systems. In particular, by comparing the RDFs, corresponding to the suggested model and obtained by implementing both the DPD and the Langevin formalisms, with the AA one at different time moments of the transient regime, we observed an almost identical behavior with insignificant differences. Thus, we demonstrated the equivalence of the two formalisms as we observed similar results for the rest of the CG models as well.

\begin{figure}[htbp]
 \centering
\includegraphics[width=0.8\textwidth]{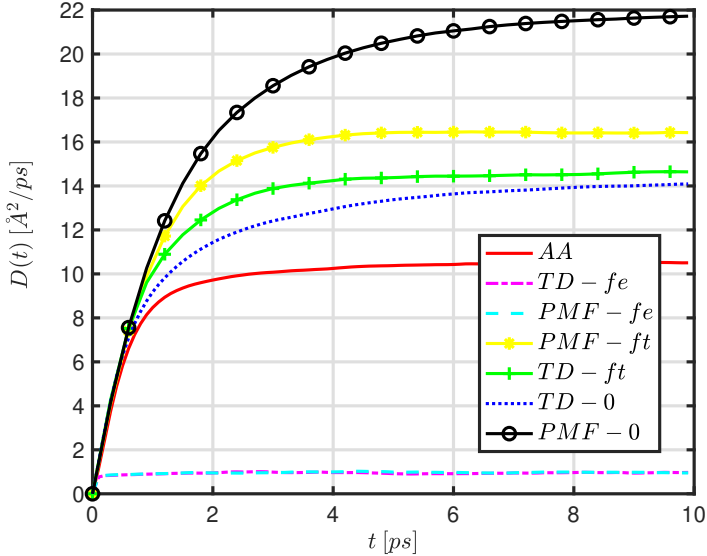}
	\caption{Diffusion coefficient of the water system, obtained from the AA and the CG simulations. \textit{PMF-0}, \textit{PMF-ft},  \textit{TD-ft}, and \textit{TD-0} models exhibit faster dynamics than \textit{TD-fe} and \textit{PMF-fe}. The DCs of \textit{TD-ft}, and \textit{TD-0} are more representative approximations of the AA DC.}
		\label{compare vacf and dc}
\end{figure}

Additionally, we obtain  information about the dynamical features of the water system by inspecting the DC. The DCs of the CG particles for all CG models as well as the AA one (analyzed at the CG level of description) are obtained from integrating the velocity autocorrelation functions of the CG particles and  are depicted in Fig. \ref{compare vacf and dc}.
We observe that the two most effective CG models among all studied in representing the DC of the system are \textit{TD-ft} and \textit{TD-0}. However, a significant deviation from the AA DC remains. The fact that the CG models \textit{PMF-ft},  \textit{TD-ft}, and \textit{TD-0} exhibit faster dynamics than the AA one indicates that the effective friction coefficient obtained using the transient AA data is smaller than should be.

\section{Discussion and Conclusions}\label{discussion-conclusions}
In this work, we presented a model for describing the transient dynamics of CG molecular systems. The model is based on a time-dependent force field to capture the memory effect of the friction kernel while retaining Markovian dynamics. For this, we introduce a Langevin-type equation in which (a) the force field explicitly depends on time and (b) the friction coefficient is time-dependent through the time interval of observation. To learn  the model parameters  we followed purely data-driven approaches: for 
the time-dependent potential' parameters  we use the PSFM approach and for the friction coefficient appropriate correlation functions. 

We applied the proposed model to two different stochastic systems. The first refers to a benchmarking example of a single particle described via the GLE.
We demonstrate that the stochastic dynamics with a time-dependent potential provide a good approximation of the short-time transient dynamics regime of the single particle moving in a box. We also examined  the effectiveness of the proposed model in terms of the correlation time of the system. For the tested values of the correlation time of the GLE, the results indicate that if the system is less correlated, the proposed  CG method is more effective in representing the CG dynamics. It can even capture the CG dynamics at a longer time interval than the transient regime. 

The second model concerns a realistic atomistic molecular model of water. The water molecules are initially in an artificial FCC configuration of low density and we simulated their evolution toward forming a water nanodroplet, under NVT conditions. Several systematic CG models are developed using data that are appropriately drawn from the detailed all-atom MD simulations. 
We demonstrate that the proposed CG model with a time-dependent potential can better approximate the CG dynamics at the short transient regime, compared to traditional Markovian models. We observe that it is more accurate in representing both the structure and the dynamics, quantified via the radial pair distribution function and the self-diffusion coefficient respectively.

At the same time, the fact that time-dependent CG models, with and without friction, exhibit  a similar efficiency emphasizes the need to study the friction kernel's effect in-depth. We also examined the effectiveness of the studied models at long time regimes, when the system reaches the equilibrium nanodroplet state. Not surprisingly, none of the transient-parametrized CG models is a representative approximation of the nanodroplet RDF at equilibrium, due to the strong fluctuations of the nanodroplet. Thus, as a future investigation, we suggest estimating the friction as an explicitly time-dependent function and incorporating it into the CG simulation to improve the CG dynamics' estimation in the short time and the long time regime as well.




A last part of the discussion, it is important to emphasize herein that, as expected, the proposed time-dependent potential we develop is specific to the particular transient experiment that we study. 
Its initial conditions are essentially a parameter of the model. Thus, the effective time-dependent potential is not transferable to systems with different initial states. However, the proposed methodology can be used to derive time-dependent effective CG models for transient systems under arbitrary initial conditions. Considering that transferability remains one of the biggest challenges in bottom-up CG models, the present paper serves as a starting point to validate the developed formalism on a realistic molecular (here water) system. 
Moreover, the specific model is transferable to different system sizes, under the same initial conditions. In other words, the derived CG potential could be used to describe the transient dynamics of a much larger water system, in a computationally efficient way, compared to the more demanding AA models. 
In future work, we also aim to examine in detail whether the time-dependent potential can be extended in time and test its validity at a longer time interval of observation than the one used to learn it, as well as to apply in non-equilibrium macromolecular systems, the atomistic simulations of which are quite often prohibitive.

Furthermore, in the context of improving the transferability of the suggested time-dependent potential, it would be of particular interest to study the inclusion of density-dependence in the approximation of the transient dynamics. In the future, together with the density-dependence, we will also further explore the thermodynamic consistency issues of time-dependent potentials.

In the next step, we also aim to extend this methodology developed for estimating the short-time transient CG dynamics to systems under non-equilibrium conditions.
In such a case  the CG force field is not necessarily a gradient and its approximation is challenging. 


\section*{Acknowledgments}

This research was supported by the Hellenic Foundation for Research and Innovation (HFRI) and the General Secretariat for Research and Technology (GSRT), under grant agreement No. 52. This work was co-funded by the European Union’s Horizon 2020 research and innovation program under grant agreement No. 810660. The research of G.B. was also partially supported by the HFRI under the HFRI Ph.D. scholarship grant agreement No. 6234.
\bibliographystyle{siamplain}
\bibliography{reference, reference1} 

\end{document}